\documentclass[11pt]{article}
\usepackage{amsfonts,amsmath,amssymb,amsthm, amscd, subcaption}
\usepackage{graphicx}
\usepackage[stable]{footmisc}
\numberwithin{equation}{section}

\usepackage{tikz}
\usepackage{tikz-cd}
\usetikzlibrary{decorations.markings}
\usetikzlibrary{arrows}
\usetikzlibrary{calc}

\newtheorem{theorem}{Theorem}[section]
\newtheorem{prop}[theorem]{Proposition}
\newtheorem{lemma}[theorem]{Lemma}
\newtheorem{cor}[theorem]{Corollary}

\theoremstyle{definition}
\newtheorem{definition}[theorem]{Definition}
\newtheorem{example}[theorem]{Example}
\newtheorem{remark}[theorem]{Remark}

\usepackage{float}

\def\({\langle \hskip -.1cm \langle} 
\def\){\rangle \hskip -.1cm \rangle} 
\def\<{{\langle}}
\def\>{{\rangle}}
\def\G{{\Gamma}}

\def\b{{\beta}}
\def\g{{\gamma}}
\def\r{{\rho}}

\def\Z{\mathbb Z}

\def\R{\mathbb R}
\def\T{\widetilde}
\def\S{{\mathbb S}}
\def\C{{\cal C}}

\def\b{\beta}

\def\t{\tau}

\def\x{\chi}
\def\G{\Gamma}

\def\s{\sigma}
\def\e{\epsilon}

\def\De{{\Delta}}
\def\Si{{\Sigma}}
\def\ti{\tilde}

\def\ni{\noindent} 

\begin{document}

\title{Core groups} 

\author{Daniel S. Silver
\and
Lorenzo Traldi
\and Susan G. Williams}

\maketitle 


\begin{abstract} The core group of a classical link was introduced independently by A.J. Kelly in 1991 and M. Wada in 1992. It  is a link invariant defined by a presentation involving the arcs and crossings of a diagram, related to Wirtinger's presentation of the fundamental group of a link complement. Two close relatives of the core group are defined by presentations involving regions rather than arcs; one of them is related to Dehn's presentation of a link group. The definitions are extended to virtual link diagrams and properties of the resulting invariants are discussed. 
\end{abstract}

\begin{center} Keywords: \textit {checkerboard coloring, classical link, Dehn presentation, virtual link};   AMS Classification: MSC: 57K10 \end{center}

\section{Introduction} \label{Intro} The core group of a classical link was introduced independently by A.J. Kelly \cite{kelly} and M. Wada 
\cite{W}. Although defined combinatorially by a presentation with generators corresponding to the arcs of a diagram of the link, it has a  topological interpretation as the free product of an infinite cyclic group and the fundamental group of the 2-fold cover of $\S^3$ branched over the link \cite{P,S,W}. Here we describe two close relatives, each group defined by a presentation obtained from a link diagram but with generators corresponding to the regions of the diagram rather than the arcs. All three groups are defined for virtual link diagrams. The first is an invariant of virtual links while the second and third are invariants for links in thickened surfaces. 

We recall some conventional notation and terminology. By a \emph{classical link} we will mean a knot or link in 3-space. (Later we will consider links in more general 3-manifolds.)  A \emph{classical link diagram} $D$ is a finite immersed collection of closed curves in the plane whose only (self-)intersections are transverse double points. An intersection is called a \emph{(classical) crossing}. At a classical crossing one of the segments is designated as the underpasser, while the other is the overpasser. They are indicated by removing a short part of the underpasser. A diagram is \emph{oriented} if a direction has been assigned to each closed component. Ambient isotopy classes of classical links correspond to equivalence classes of classical link diagrams under classical Reidemeister moves (see, for example, \cite{lickorish}).  

A \emph{virtual link diagram} is a link diagram as above but which may also contain a second type of intersection called a \emph{virtual crossing}, and indicated by a small circle enclosing it. (A virtual link diagram is not required to have any virtual crossings.) A \emph{virtual link} is an equivalence class of virtual link diagrams under generalized Reidemeister moves, a finite set of local changes of a diagram that includes Reidemeister moves. The concept is due to L. Kauffman and details can be found in \cite{kauffman}. If two classical link diagrams are equivalent under generalized Reidemeister moves, then they are equivalent under classical Reidemeister moves \cite{gpv}. Consequently, virtual link theory extends the classical theory. By a \emph{virtual link} we will mean a classical or virtual link. \bigskip

We begin with the definition of Kelly and Wada's core group, which we will call the \emph{arc core group}. If $D$ is a virtual link diagram then an \emph{arc} of $D$ is a segment that extends in both directions until it arrives at a classical crossing where it is an underpasser. The set of arcs of $D$ is denoted $A(D)$, and the set of classical crossings of $D$ is denoted $C(D)$.

\begin{definition}
\label{arccore}
Let $D$ be a virtual link diagram, and let $\{g_a \mid a \in A(D)\}$ be a set of symbols in one-to-one correspondence with $A(D)$. If $c\in C(D)$ is a classical crossing with overpassing arc $a_1(c)$ and underpassing arcs $a_2(c)$ and $a_3(c)$, let $r_c=g_{a_1(c)}g^{-1}_{a_2(c)}g_{a_1(c)}g^{-1}_{a_3(c)}$. Then the \emph{arc core group} $AC(D)$ is the group with presentation $\langle \{g_a \mid a \in A(D)\};\{r_c \mid c \in C(D)\} \rangle$.
\end{definition}

The notation $AC(D)$ is intended to emphasize that the arc core group is defined using $A(D)$. Notice that $AC(D)$ is well defined, because the group-theoretic significance of a crossing relator $r_c$ is not changed if the indices of the underpassing arcs $a_2(c),a_3(c)$ are interchanged. It is well known (and not difficult to verify) that $AC(D)$ is invariant under generalized Reidemeister moves. Hence it is a virtual link invariant $AC(L)$. 

The other members of the core group family that we discuss are defined using the \emph{regions} of a link diagram, rather than the arcs. For a classical link diagram $D$, the regions are defined to be the connected components of the complement in the plane of the union of immersed circles from which $D$ is built. For a virtual link diagram $D$, the regions are defined to be the connected components of the complement of the union of immersed circles in the \emph{abstract link diagram} standardly associated with $D$ \cite{KK}. An abstract link diagram is a compact orientable surface that has as a deformation retraction  the underlying 4-valent graph of $D$, with vertices corresponding to classical crossings (and virtual crossings ignored). 

For a classical link diagram, the classical and virtual notions of ``regions'' are not the same. There are two reasons for this disparity. 

(1.) Recall that a link diagram $D$ is \emph{split} if there is an embedded circle in the plane with nonempty parts of $D$ both inside and outside. The abstract link diagram associated with a split diagram $D$ is the disjoint union of the abstract link diagrams associated with the non-split subdiagrams whose split union is $D$. Therefore no two of these subdiagrams share any region. However, if $D$ is a split classical diagram considered in the plane, then the non-split subdiagrams that constitute $D$ are connected to each other by shared complementary regions. As a result, if $D$ is the split union of $k$ disjoint non-split subdiagrams, then $D$ will have $k-1$ fewer regions according to the classical definition than it has according to the virtual definition. 

(2.) Even for a non-split classical diagram, there is a difference between the classical regions (the bounded regions are disks) and the virtual regions (which are annuli). This difference has no significance in our discussion, as there is a natural correspondence between the classical regions and the virtual regions.

The disparity between the classical and virtual notions of ``regions'' is reflected in the following.

\begin{definition}
\label{regiondef} If $D$ is a classical link diagram then $R(D)$ denotes the set of regions of $D$ according to the classical definition, and $R_v(D)$ denotes the set of regions of $D$ according to the virtual definition. If $D$ is a non-classical virtual link diagram, then $R(D)$ denotes the set of regions of $D$ according to the virtual definition.
\end{definition}

\begin{figure} [bht]
\centering
\begin{tikzpicture}
\draw [thick] (-6,0) -- (-4.2,0);
\draw [thick] (-3.8,0) -- (-2,0);
\draw [thick] (-4,1.5) -- (-4,-1.5);
\draw [thick] (6,0) -- (4.2,0);
\draw [thick] (6,0.5) -- (4.7,0.5);
\draw [thick] (4.7,1.5) -- (4.7,0.5);
\draw [thick] (6,-0.5) -- (4.7,-0.5);
\draw [thick] (4.7,-1.5) -- (4.7,-0.5);
\draw [thick] (3.3,1.5) -- (3.3,0.5);
\draw [thick] (2,0.5) -- (3.3,0.5);
\draw [thick] (3.3,-1.5) -- (3.3,-0.5);
\draw [thick] (2,-0.5) -- (3.3,-0.5);
\draw [thick] (3.8,0) -- (2,0);
\draw [thick] (4,1.5) -- (4,-1.5);
\node at (-3.6,0.35) {$W$};
\node at (-4.4,0.35) {$V$};
\node at (-4.4,-0.35) {$Y$};
\node at (-3.6,-0.35) {$X$};
\node at (3.6,0.35) {$V$};
\node at (4.4,0.35) {$W$};
\node at (4.4,-0.35) {$X$};
\node at (3.6,-0.35) {$Y$};
\end{tikzpicture}
\caption{The regions incident at a classical crossing in a classical diagram (on the left) or an abstract link diagram (on the right).}
\label{regionfig}
\end{figure}
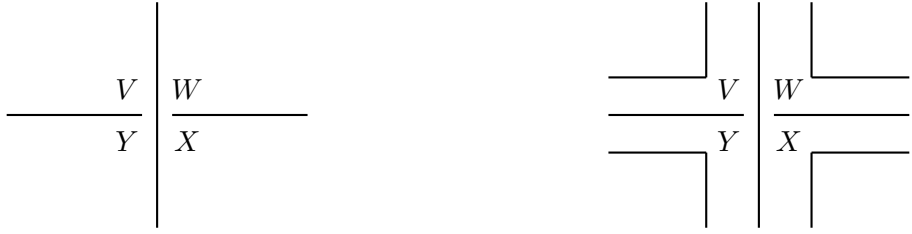

\begin{definition}
\label{regioncore} Let $D$ be a virtual link diagram, let $R(D)$ be the set of regions of $D$, and let $\{\g_R \mid R \in R(D)\}$ be a set of symbols in one-to-one correspondence with $R(D)$. If $c\in C(D)$ is a classical crossing as pictured in Figure \ref{regionfig}, let $\r_c=\g_V^{}\g^{-1}_W\g_Y^{}\g^{-1}_X$. Then the \emph{first regional core group} $RC(D)$ is the group with presentation $\langle \{\g_R \mid R \in R(D)\};\{\r_c \mid c \in C(D)\} \rangle$.
\end{definition}

The significance of a relator $r$ in a group presentation is unchanged if $r$ is replaced by its inverse, or if the terms in $r$ or $r^{-1}$ are cyclically permuted. It follows that if $c$ is the crossing pictured in Figure \ref{regionfig}, then the significance of $\r_c$ is unchanged if the regions are relabeled as follows: $V$ is any one of the four regions; $W$ is the region across the overpassing arc from $V$; $X$ is the region diagonally across the crossing from $V$; and $Y$ is the region across an underpassing arc from $V$. 

In Section \ref{proof1} we prove:

\begin{theorem}
\label{main1}
If $D$ is a classical link diagram, then $RC(D)$ is isomorphic to the free product $ \Z * AC(D)$. 
\end{theorem}

Our second regional core group is defined using the following notion of crossing index. 

\begin{definition}
\label{stindex}
If $c$ is a classical crossing as illustrated in Figure \ref{regionfig}, then the \emph{crossing indices} of $c$ are $\eta_c(V,X)=\eta_c(X,V)=-1$ and $\eta_c(W,Y)=\eta_c(Y,W)=1$.
\end{definition}

The reader familiar with the account of the classical Goeritz matrix in Lickorish's textbook \cite{lickorish} will recognize that for each pair of diagonally opposite regions, the indices $\eta_c$ agree with the Goeritz index of $c$ when those two regions are shaded.

\begin{definition}
\label{rregioncore}
Let $D$ be a virtual link diagram, and let $\{g_R \mid R \in R(D)\}$ be a set of symbols in one-to-one correspondence with $R(D)$. For each region $R \in R(D)$, let $c_1, \dots,c_k$ be the classical crossings of $D$ incident on $R$, listed in the order we encounter them if we walk counterclockwise around the boundary of $R$; for $i \in \{1, \dots,k\}$ let $Q_i$ be the region diagonally across from $R$ at $c_i$. Then the \emph{second regional core group} $RRC(D)$ is the group with presentation $\langle \{g_R \mid R \in R(D)\};\{r_R \mid R \in R(D)\} \rangle$, where the relator $r_R$ is given by 
\[
r_R = \prod_{i=1}^k (g^{-1}_Rg_{Q_i})^{\eta_{c_i}(R,Q_i)}.
\]
\end{definition}

\begin{theorem}
\label{main2}
Let $D$ be a classical link diagram that is the split union of $k$ disjoint non-split diagrams, and let $F_{k+1}$ be a free group of rank $k+1$. Then the free products $ RC(D) * RC(D)$ and $F_{k+1}*RRC(D)$ are isomorphic.
\end{theorem}

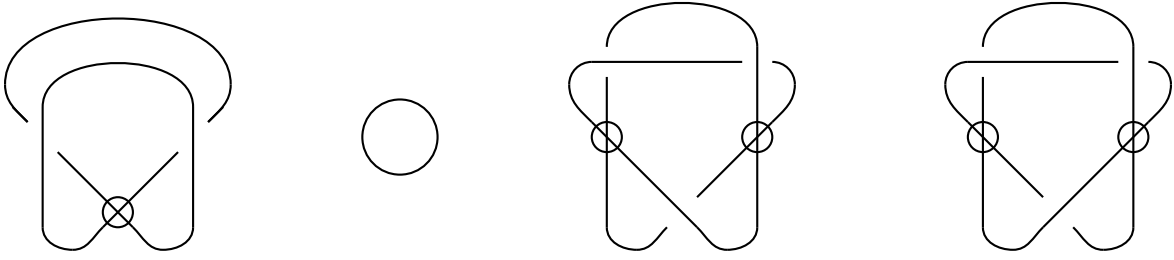
\begin{figure} [bth]
\centering
\begin{tikzpicture} [>=angle 90]
\draw [thick] (11.2,1.7) to [out=45, in=-90] (11.5,2.2);
\draw [thick] (11.2,2.5) to [out=0, in=90] (11.5,2.2);
\draw [thick] (11.2,1.7) --(9.8,.3);
\draw [thick] (8.8,1.7) --(9.8,.7);
\draw [thick] (10.2,.3) to [out=-45, in=180] (10.6,0);
\draw [thick] (11,.3) to [out=-90, in=0] (10.6,0);
\draw [thick] (10.8,2.5) --(8.8,2.5);
\draw [thick] (8.8,1.7) to [out=135, in=-90] (8.5,2.2);
\draw [thick] (8.8,2.5) to [out=180, in=90] (8.5,2.2);
\draw [thick] (11,2.7) to [out=90, in=90] (9,2.7);
\draw [thick] (11,2.7) -- (11,.3);
\draw [thick] (9,2.3) -- (9,.3);
\draw [thick] (9.8,.3) to [out=-135, in=0] (9.4,0);
\draw [thick] (9,.3) to [out=-90, in=180] (9.4,0);
\draw [thick] (6.2,1.7) to [out=45, in=-90] (6.5,2.2);
\draw [thick] (6.2,2.5) to [out=0, in=90] (6.5,2.2);
\draw [thick] (6.2,1.7) --(5.2,.7);
\draw [thick] (3.8,1.7) --(5.2,.3);
\draw [thick] (5.2,.3) to [out=-45, in=180] (5.6,0);
\draw [thick] (6,.3) to [out=-90, in=0] (5.6,0);
\draw [thick] (5.8,2.5) --(3.8,2.5);
\draw [thick] (3.8,1.7) to [out=135, in=-90] (3.5,2.2);
\draw [thick] (3.8,2.5) to [out=180, in=90] (3.5,2.2);
\draw [thick] (6,2.7) to [out=90, in=90] (4,2.7);
\draw [thick] (6,2.7) -- (6,.3);
\draw [thick] (4,2.3) -- (4,.3);
\draw [thick] (4.8,.3) to [out=-135, in=0] (4.4,0);
\draw [thick] (4,.3) to [out=-90, in=180] (4.4,0);
\draw[fill=none] [thick](4,1.5) circle (0.2);
\draw[fill=none] [thick](6,1.5) circle (0.2);
\draw[fill=none] [thick](9,1.5) circle (0.2);
\draw[fill=none] [thick](11,1.5) circle (0.2);
\draw[fill=none] [thick](1.25,1.5) circle (.5);
\draw [thick] (11.2-12.5,1.7) to [out=45, in=-90] (11.5-12.5,2.2);
\draw [thick] (10.8-12.5,1.3) --(9.8-12.5,.3);
\draw [thick] (11.2-12.5,1.7) --(11.4-12.5,1.9);
\draw [thick] (9.2-12.5,1.3) --(10.2-12.5,.3);
\draw [thick] (8.8-12.5,1.7) --(8.6-12.5,1.9);

\draw [thick] (10.2-12.5,.3) to [out=-45, in=180] (10.6-12.5,0);
\draw [thick] (11-12.5,.3) to [out=-90, in=0] (10.6-12.5,0);

\draw [thick] (8.8-12.5,1.7) to [out=135, in=-90] (8.5-12.5,2.2);
\draw [thick] (11.5-12.5,2.2) to [out=90, in=90] (8.5-12.5,2.2);

\draw [thick] (11-12.5,1.9) -- (11-12.5,.3);
\draw [thick] (9-12.5,1.9) -- (9-12.5,.3);
\draw [thick] (9-12.5,1.9) to [out=90, in = 90] (11-12.5,1.9);
\draw [thick] (9.8-12.5,.3) to [out=-135, in=0] (9.4-12.5,0);
\draw [thick] (9-12.5,.3) to [out=-90, in=180] (9.4-12.5,0);
\draw[fill=none] [thick](-2.5,.5) circle (0.2);
\end{tikzpicture}
\caption{The group triples $(AC,RC,RRC)$ for these four diagrams are $(\mathbb Z, \mathbb Z,\mathbb Z*\mathbb Z)$, $(\mathbb Z, \mathbb Z*\mathbb Z,\mathbb Z*\mathbb Z)$, $(\mathbb Z, \mathbb Z*\mathbb Z,\mathbb Z*\mathbb Z*\mathbb Z_2)$ and $(\mathbb Z * \mathbb Z_3, \mathbb Z*\mathbb Z,\mathbb Z*\mathbb Z*\mathbb Z_2)$, from left to right.}
\label{virttref}
\end{figure}

Theorems \ref{main1} and \ref{main2} imply that for non-split classical link diagrams, any one of the three groups $AC(D), RC(D), RRC(D)$ determines the other two. The situation for virtual diagrams is quite different: in general no two of the groups determine the third. See Figure  \ref{virttref} for examples.

Recall that two virtual link diagrams are equivalent if one can be transformed into the other through some finite sequence of detour moves (which involve only virtual crossings) and Reidemeister moves $\Omega.1, \Omega.2, \Omega.3$ (which involve only classical crossings). The groups $AC(D), RC(D)$ and $RRC(D)$ are all invariant under $\Omega.1$ and $\Omega.3$ moves, and also under detour moves. For classical diagrams, $AC(D)$ and $RC(D)$ are invariant under $\Omega.2$ moves, and if $D$ is the split union of $k$ non-split diagrams, then $F_{k+1}*RRC(D)$ is invariant under $\Omega.2$ moves. 

For virtual link diagrams $D$, the group $AC(D)$ remains invariant under $\Omega.2$ moves, but the regional groups are not invariant. Nevertheless, $RC(D)$ and $RRC(D)$ provide invariants of links in thickened surfaces. See Sections \ref{Reidemeister} and \ref{Dehnc}.

\section{Proof of Theorem \ref{main1}}
\label{proof1}

We recall that a virtual link diagram $D$ is \emph{checkerboard colorable} if it is possible to shade its regions in such a way that at every non-crossing point of an arc $a$ of $D$, there is a shaded region on one side of $a$ and an unshaded region on the other side of $a$. All classical diagrams are checkerboard colorable.

\begin{prop}
\label{themapf}
If $D$ is a checkerboard colorable virtual link diagram, then there is a function $f:A(D) \to RC(D)$ defined as follows: If $a \in A(D)$ has a shaded region $S$ on one side and an unshaded region $U$ on the other side, then $f(a)=\g_U \g^{-1}_S$.
\end{prop}

\begin{proof}
Well-definedness of $f$ requires that if $c$ is a crossing with  overpassing arc $a$, then we have the same value for $\g_U^{} \g^{-1}_S$ on both sides of $c$. 

Suppose that $a$ is the overpassing arc at the crossing $c$ pictured in Figure \ref{regionfig}.
If the shaded regions are $W$ and $Y$ then $\g_U^{} \g^{-1}_S=\g_V^{}  \g^{-1}_W$ on one side of $c$, and $\g_U^{}  \g^{-1}_S=\g_X^{}  \g^{-1}_Y$ on the other side of $c$. The crossing relator $\g_V^{}  \g^{-1}_W\g_Y^{}  \g^{-1}_X$ tells us that $\g_V^{}  \g^{-1}_W=\g_X^{}  \g^{-1}_Y$.
Similarly, if the shaded regions are $V$ and $X$ then $\g_U^{}  \g^{-1}_S=\g_W^{}  \g^{-1}_V$ on one side of $c$, and $\g_U^{}  \g^{-1}_S=\g_Y^{} \g^{-1}_X$ on the other side of $c$. The crossing relator $\g_V^{}  \g^{-1}_W\g_Y^{}  \g^{-1}_X$ tells us that $\g_W^{}  \g^{-1}_V=\g_Y^{}  \g^{-1}_X$. \end{proof}

\begin{prop}
\label{the homf}
If $D$ is a checkerboard colorable virtual link diagram then there is a homomorphism $AC(D) \to RC(D)$ given by $g_a \mapsto f(a) \thickspace \allowbreak \forall a \in A(D)$. 
\end{prop}
\begin{proof}
We must show that the crossing relators of Definition \ref{arccore} are mapped to $1$ in $RC(D)$.

Consider a crossing $c$ as pictured in Figure \ref{regionfig}, with overpasser $a$ and underpassers $b$ (on the left) and $b'$ (on the right). Then $r_c=g_ag^{-1}_bg_ag_{b'}^{-1}$. If the shaded regions are $V$ and $X$, the image of $r_c$ is
\[
f(g_a)f(g_b)^{-1}f(g_a)f(g_{b'})^{-1}=\g_W^{} \g^{-1}_V(\g_Y^{} \g^{-1}_V)^{-1}\g_Y^{} \g^{-1}_X(\g_W^{} \g^{-1}_X)^{-1}=1\text{,}
\]
as required. Similarly, if the shaded regions are $W$ and $Y$ then the image of $r_c$ is
\[
f(g_a)f(g_b)^{-1}f(g_a)f(g_{b'})^{-1}=\g_X^{} \g^{-1}_Y(\g_V^{} \g^{-1}_Y)^{-1}\g_V^{} \g^{-1}_W(\g_X^{} \g^{-1}_W)^{-1}=1.
\]
\end{proof}
We abuse notation slightly by using the letter $f$ to denote the homomorphism $AC(D) \to RC(D)$ given by $f(g_a)=f(a) \thickspace \allowbreak \forall a \in A(D)$.

\begin{lemma}
\label{grouplem}
Let $G$ be a group given by a presentation $\langle S; R \rangle$ in which every relator $r \in R$ is of the form $s_1s_2^{-1} \dots s_{2k-1}s_{2k}^{-1}$ for some $k \in \mathbb N$ and some $s_1,\dots, s_{2k} \in S$. Let $s_0$ be any fixed element of $S$, let $( s_0 )$ be the cyclic subgroup of $G$ generated by $s_0$, and let $H$ be the subgroup of $G$ generated by $\{ss^{-1}_0 \mid s \neq s_0 \in S\}$. Then $( s_0 )$ is infinite, and $G$ is the internal free product $( s_0 ) * H$.
\end{lemma}
\begin{proof}
For each $s \in S$, let $s'=ss^{-1}_0$. If $S'=\{s_0\} \cup \{s' \mid s \neq s_0 \in S\}$ then the group presentation $\langle S; R \rangle$ can be rewritten as $\langle S'; R' \rangle$, where for each relator $r=s_1s_2^{-1} \dots s_{2k-1}s_{2k}^{-1} \in R$, $r'$ is the word obtained from $s'_1(s'_2)^{-1} \dots s'_{2k-1}(s'_{2k})^{-1}$ by removing all occurrences of $s'_0$. The lemma follows, because $s_0$ does not appear in any relator $r' \in R'$. \end{proof}

\begin{cor}
\label{interchange}
In the situation of Lemma \ref{grouplem}, suppose $s_0,s_1 \in S$. Then (a) $\{ss^{-1}_0 \mid s \neq s_0 \in S\}$ and $\{ss^{-1}_1 \mid s \neq s_1 \in S\}$ generate the same subgroup $H$ of $G$, and (b) $G$ has an automorphism $\varphi$ with $\varphi(s_0)=s_1$ and $\varphi(h)=h \thickspace \allowbreak \forall h \in H$.
\end{cor}
\begin{proof} For (a), observe that for every $s \in S$, $ss_0^{-1} = (ss_1^{-1})(s_0s_1^{-1})^{-1}$ and $ss_1^{-1} = (ss_0^{-1})(s_1s_0^{-1})^{-1}$. For (b), 
Lemma \ref{grouplem} tells us that 
$( s_0 )$ and $( s_1 )$ are both infinite cyclic groups, so there is an isomorphism $\phi:( s_0 ) \to ( s_1 )$ with $\phi(s_0) = s_1$. Lemma \ref{grouplem} also tells us that $( s_0 ) * H =G=( s_1 ) * H$. The automorphism $\varphi$ of $G$ is obtained by combining $\phi$ with the identity map of $H$.
\end{proof}
Notice that Lemma \ref{grouplem} and Corollary \ref{interchange} apply to both $AC(D)$ and $RC(D)$.

\begin{lemma}
\label{flem}
Let $D$ be a classical link diagram or a non-split checkerboard colorable virtual link diagram, and let $R_0$ be any fixed region of $D$. Then the subgroup $H$ of $RC(D)$ generated by $\{ \g_R^{}  \g_{R_0}^{-1} \mid R \neq R_0 \in R(D) \}$ is the same as $f(AC(D))$.
\end{lemma}
\begin{proof}
For any arc $a \in A(D)$, there are regions $U,S \in R(D)$ with 
\[
f(g_a)=\g_U^{}  \g ^{-1}_S=(\g_U^{}  \g^{-1}_{R_0})(\g_S^{}  \g^{-1}_{R_0})^{-1} \in H.
\]
It follows that $f(AC(D)) \subseteq H$.

For the opposite inclusion, let $R$ be any region of $D$. The hypothesis that $D$ is non-split or classical implies that there is a sequence $R_0,R_1, \dots, R_k=R$ of regions of $D$ such that for each $i \in \{1, \dots, k \}$, there is an arc $a_i\in A(D)$ with $R_{i-1}$ on one side and $R_i$ on the other side. It follows that $\g_{R_i} \g^{-1}_{R_{i-1}} = f(g_{a_i})^{\pm 1} \in f(AC(D))$ for each $i$, and hence 
\[
\g_R \g^{-1}_{R_0} = (\g_{R_k} \g^{-1}_{R_{k-1}})( \g_{R_{k-1}} \g^{-1}_{R_{k-2}})\dots (\g_{R_1} \g^{-1}_{R_0}) \in f(AC(D)).
\]
It follows that $f(AC(D)) \supseteq H$. \end{proof}

Combining Lemmas \ref{grouplem} and \ref{flem}, we deduce the following.

\begin{prop}
\label{freeproduct}
If $D$ is a classical link diagram or a non-split checkerboard colorable virtual link diagram, and $R_0$ is any region of $D$, then $( \g_{R_0} )$ is an infinite cyclic subgroup of $RC(D)$, and $RC(D)$ is the internal free product $( \g_{R_0} ) * f(AC(D))$.
\end{prop}

To complete the proof of Theorem \ref{main1} we verify that if $D$ is a classical link diagram, then the homomorphism $f:AC(D) \to RC(D)$ is injective. We do this by proving the following.

\begin{theorem}
\label{linverse}
If $D$ is a classical link diagram, then there is a homomorphism $h:RC(D) \to AC(D)$ with $hf(g_a)=g_a \allowbreak \thickspace \forall a \in A(D)$.
\end{theorem}

\begin{proof}
We begin the process of defining $h$ by choosing a fixed shaded region of $D$, $S_0$. 

If $S$ is any shaded region of $D$, then we can find a path $P(S)$ from $S_0$ to $S$ in the plane, which intersects $D$ only finitely many times, and does not pass through any crossing. Following the path $P(S)$ from $S_0$ to $S$, we obtain a sequence $S_0,a_0,U_0,b_0,S_1,a_1,U_1, \dots, U_k,b_k,S_{k+1}=S$ such that every $S_i$ is a shaded region of $D$, every $U_i$ is an unshaded region of $D$, every $a_i$ is an arc of $D$ with $S_i$ on one side and $U_i$ on the other side, and every $b_i$ is an arc of $D$ with $U_i$ on one side and $S_{i+1}$ on the other side. Denote this sequence $WP(S)$, and let $h(WP(S))=g^{-1}_{b_k}g_{a_k}g^{-1}_{b_{k-1}} \dots g^{-1}_{b_0}g_{a_0}\in AC(D)$. 

We claim that if $P'(S)$ is some other path from $S_0$ to $S$, then $h(WP(S))=h(WP'(S))$. To verify this claim, consider that the fact that $\R^2$ is simply connected implies that $P$ can be continuously deformed into $P'$. During such a continuous deformation, the sequence $WP$ changes only when the deformation passes through a crossing. 

To verify the claim, suppose we label the arcs incident at a crossing $c$ as $a,a',a'',a'''$ in either clockwise or counterclockwise order, with any arc serving as $a$. (The two labels corresponding to the overpasser will represent the same arc, so $a=a''$ or $a'=a'''$.) Then the crossing relation $1=g_{a_1(c)}g^{-1}_{a_2(c)}g_{a_1(c)}g^{-1}_{a_3(c)}$ of Definition \ref{arccore} will imply that the equalities 
\[
g_{a} = g_{a'}g^{-1}_{a''}g_{a'''} \qquad \text{   and   } \qquad g^{-1}_{a} = g^{-1}_{a'}g_{a''}g^{-1}_{a'''}
\]
both hold. If $a$ appears in $WP(S)$ as one of the $a_i$, then the effect on $h(WP(S))$ when the deformation passes through $c$ is to replace $g_a$ with $g_{a'}g^{-1}_{a''}g_{a'''}$, and perhaps to cancel some terms. If $a$ appears in $WP(S)$ as one of the $b_i$, then the effect on $h(WP(S))$ when the deformation passes through $c$ is to replace $g^{-1}_a$ with $g^{-1}_{a'}g_{a''}g^{-1}_{a'''}$, and perhaps to cancel some terms. Either way, the value of $h(WP(S))$ is not changed.

Similarly, if $U$ is an unshaded region of $D$ then we can find a path $P(U)$ from $S_0$ to $U$ in the plane, which intersects $D$ only finitely many times, and does not pass through any crossing. Following the path $P(U)$ from $S_0$ to $U$ provides a sequence $S_0,a_0,U_0,b_0,S_1,a_1,U_1, \dots, a_k, U_k=U$ such that every $S_i$ is a shaded region of $D$, every $U_i$ is an unshaded region of $D$, every $a_i$ is an arc of $D$ with $S_i$ on one side and $U_i$ on the other side, and every $b_i$ is an arc of $D$ with $U_i$ on one side and $S_{i+1}$ on the other side. Denote this sequence $WP(U)$, and let $h(WP(U))=g_{a_k}g^{-1}_{b_{k-1}} \dots g^{-1}_{b_0}g_{a_0}\in AC(D)$. The same kind of argument given above shows that $h(WP(U))$ does not depend on the choice of the path $P(U)$.

We now have a well-defined function $h:R(D) \to AC(D)$. To verify that $h$ defines a homomorphism $RC(D) \to AC(D)$, suppose $c$ is a crossing as pictured in Figure \ref{regionfig}. In the introduction we mentioned that the choice of region labels $V,W,X,Y$ in Figure \ref{regionfig} is quite flexible; $V$ can be any one of the four regions, so long as $W$ neighbors $V$ across the overpassing arc and $Y$ neighbors $V$ across an underpassing arc. In particular, we may choose the region labels so that $V$ is unshaded. Let $a$ be the overpassing arc of $c$. We may assume that the path $P(W)$ is obtained by first following $P(V)$, and then crossing $a$ from $V$ to $W$. Then 
\[
h(WP(V))h(WP(W))^{-1}=(g_{a_k}g^{-1}_{b_{k-1}} \dots g^{-1}_{b_0}g_{a_0})(g^{-1}_ag_{a_k}g^{-1}_{b_{k-1}} \dots g^{-1}_{b_0}g_{a_0})^{-1}=g_a.
\]
A similar calculation shows that $h(WP(X))h(WP(Y))^{-1}=g_a$ too. Then 
\[
h(WP(V))h(WP(W))^{-1}h(WP(Y))h(WP(X))^{-1}=g_ag^{-1}_a=1 \text{,}
\]
so the image under $h$ of the crossing relator $\r_c=\g_V^{} \g^{-1}_W \g_Y\g^{-1}_X$ is trivial. It follows that $\g_R \mapsto h(WP(R)$ defines a homomorphism $RC(D) \to AC(D)$, which we also denote $h$.

To complete the proof, we need to show that if $a \in A(D)$ then $h(f(g_a))=g_a$. The argument is the same as the calculation of $h(WP(V))h(WP(W))^{-1}=g_a$ in the preceding paragraph. \end{proof}

\section{Proof of Theorem \ref{main2}}
\label{proof2}
Before beginning the proof, observe that if $D$ is a checkerboard colorable link diagram, then for each region $R \in R(D)$, the relator $r_R$ of Definition \ref{rregioncore} involves only generators corresponding to regions with the same shaded/unshaded status as $R$ itself. It follows that the group $RRC(D)$ is the internal free product $RRC(D)=RRC_S(D) * RRC_U(D)$, where the subgroup $RRC_S(D)$ corresponds to the shaded regions and the subgroup $RRC_U(D)$ corresponds to the unshaded regions.

Now, let $D$ be a classical link diagram.  There is a \emph{shaded checkerboard graph} $\Gamma_s(D)$, with a vertex for each shaded region and an edge for each crossing; the edge corresponding to the crossing $c$ is incident on the vertex or vertices corresponding to shaded region(s) incident at $c$. There is also an \emph{unshaded checkerboard graph} $\Gamma_u(D)$, with vertices for the unshaded regions of $D$ and edges for crossings of $D$, defined in the same way. We use $\b_s(D)$ to denote the number of connected components of $\G_s(D)$, and $\b_u(D)$ to denote the number of connected components of $\G_u(D)$.

We can eliminate most of the generators of $RC(D)$ corresponding to shaded regions, as follows. Choose shaded regions $S_1,\dots, S_{\b_s(D)}$, one in each connected component of the graph $\G_s(D)$. Suppose $S$ is a shaded region that is in the same connected component of $\G_s(D)$ as $S_n$, but is not equal to $S_n$. Then there is a path $P(S)$ from $S_n$ to $S$ in the plane, which avoids the unshaded regions of $D$. Say $P(S)$ is $S_n=\Sigma_0,c_1,\Sigma_1, \dots, c_k,\Sigma_k=S$, where each $\Sigma_i$ is a shaded region and each $c_i$ is a crossing. For $1\leq i \leq k$ let $U_i$ be the unshaded region on the left of $c_i$ (as seen by an observer standing in $\Sigma_i$ and facing $\Sigma_{i+1}$), let $V_i$ be the unshaded region on the right, and let $\eta_i=\eta_{c_i}(U_i,V_i)$ according to Definition \ref{stindex}. Also, let 
\[
\g(P(S))=\prod_{i=1}^k (\g^{-1}_{U_{i}}\g_{V_{i}})^{\eta_i}.
\]

According to Definition \ref{regioncore}, $\g_{\Sigma_{i}}\g^{-1}_{U_{i}}\g_{V_{i}}\g^{-1}_{\Sigma_{i+1}}=1$ if $\eta_i=1$, and $\g_{U_{i}}\g^{-1}_{\Sigma_{i+1}}\g_{\Sigma_{i}}\g^{-1}_{V_{i}}=1$ if $\eta_i=-1$. Either way, 
\[
\g^{-1}_{\Sigma_{i}}\g_{\Sigma_{i+1}}=(\g^{-1}_{U_{i}}\g_{V_{i}})^{\eta_i}.
\]
It follows that
\[
\g_S = \g_{\Sigma_{0}}\g^{-1}_{\Sigma_{0}}\g_{\Sigma_{1}}\g^{-1}_{\Sigma_{1}}\g_{\Sigma_{2}} \dots \g^{-1}_{\Sigma_{k-1}}\g_{\Sigma_{k}}=\g_{S_n} \cdot \prod_{i=1}^k (\g^{-1}_{U_{i}}\g_{V_{i}})^{\eta_i}=\g_{S_n} \cdot \g(P(S)).
\]

If we use these relations $\g_S = \g_{S_n} \cdot \g(P(S))$ to eliminate every generator $\g_S$ with $S \notin \{S_1, \dots, \allowbreak S_{\b_s(D))}\}$, then we will be left with a presentation of $RC(D)$ that has generators from the set $\{ \g_U \mid U \text{ is an unshaded region of }D\} \cup  \{\g_{S_1}, \dots, \g_{S_{\b_s(D))}}\}$, and has relations stating that whenever $P(S)$ and $P'(S)$ are paths from an $S_n$ to the same unshaded region $S$, $\g(P(S)) = \g(P'(S))$. Notice that the generators $\g_{S_n}$ do not appear in any of these relations.

The relations $\g(P(S))=\g(P'(S))$ may be equivalently described this way: $\g(P)=1$ whenever $P$ is a closed path in $\G_s(D)$ that begins and ends in $S_n$. It is well known that the closed paths in the planar graph $\G_s(D)$ are generated by the boundary circuits of the unshaded regions. That is, to generate the relations $\g(P)=1$ it suffices to consider closed paths of the form $P=P_1P_2P_3$, where $P_2$ is the boundary of an unshaded region $U$ and $P_3$ is the same path as $P_1$, traversed in the opposite direction. Notice that this relationship between $P_1$ and $P_3$ implies that when a crossing is traversed as part of $P_3$, $\eta_i$ remains the same but the regions denoted $U_i$ and $V_i$ are interchanged. It follows that 
\[
\g(P)=\g(P_1)\g(P_2)\g(P_3) = \g(P_1)\g(P_2)\g(P_1)^{-1}
\]
so as a relator, $\g(P)$ is equivalent to $\g(P_2)$, which is the same as the relator $r_U$ of Definition \ref{rregioncore}.

Putting all of this together, we obtain an isomorphism
\[
RC(D) \cong  RRC_U(D) * F_{\b_s(D)}.
\]
Reversing the roles of the shaded and unshaded regions, the same argument leads to an isomorphism
\[
RC(D) \cong  RRC_S(D) * F_{\b_u(D)}.
\]
Combining the two isomorphisms, we have 
\[
RC(D)*RC(D) \cong  RRC_U(D)*RRC_S(D) * F_{\b_s(D)+\b_u(D)}.
\]
Theorem \ref{main2} follows, because $RRC_U(D)*RRC_S(D)\cong RRC(D)$ and if $D$ is the split union of $k$ non-split diagrams, then $\b_s(D)+\b_u(D)=k+1$.  

\begin{remark}In \cite{STW} it is shown that if $L$ is a classical link with Goeritz matrix $G$, and $\bf 0$ denotes a zero matrix with $\b_s(D)$ columns, then there is a rewritten version of the Dehn presentation of $\pi_1(\S^3-L)$ that yields the matrix $\begin{pmatrix} G & \bf 0 \end{pmatrix}$ through Fox's free differential calculus. Here the algebra is not so complicated, as the free differential calculus is not required: the natural presentation of $ RRC_U(D) * F_{\b_s(D)}$ yields $\begin{pmatrix} G & \bf 0 \end{pmatrix}$ through simple abelianization.\end{remark}

\section{Reidemeister moves}
\label{Reidemeister}

We leave the rather mechanical proofs of Propositions \ref{R1} and \ref{R2} to the reader. 

\begin{prop} \label{R1} For any virtual link diagram $D$, the arc core group $AC(D)$ is invariant under detour moves and Reidemeister moves. That is, $AC(D)$ is an invariant of virtual link type.
\end{prop}

\begin{prop} \label{R2} For any virtual link diagram $D$, the region core groups $RC(D)$ and $RRC(D)$ are invariant under detour moves and Reidemeister moves of types $\Omega.1$ and $\Omega.3$. 
\end{prop}

\begin{figure} [bth]
\centering
\begin{tikzpicture} [>=angle 90]
\draw [thick] (11.2,1.7) to [out=45, in=-90] (11.5,2.2);
\draw [thick] (11.2,2.5) to [out=0, in=90] (11.5,2.2);
\draw [thick] (8.8,1.7) --(10.2,.3);
\draw [thick] (10.2,.3) to [out=-45, in=180] (10.6,0);
\draw [thick] (11,.3) to [out=-90, in=0] (10.6,0);
\draw [thick] (11.2,2.5) --(8.8,2.5);
\draw [thick] (8.8,1.7) to [out=135, in=-90] (8.5,2.2);
\draw [thick] (8.8,2.5) to [out=180, in=90] (8.5,2.2);
\draw [thick] (11,2.7) to [out=90, in=90] (9,2.7);
\draw [thick] (9,2.7) -- (9,2.5);
\draw [thick] (11,2.3) -- (11,.3);
\draw[fill=none] [thick](9,2.5) circle (0.2);
\draw[fill=none] [thick](11,1.5) circle (0.2);
\draw [thick] (11.2,1.7)--(10.7,1.2);
\draw [thick] (9,2.5) to [out=-90, in = -135] (10.7,1.2);
\draw [thick] (11.2-4.2,1.7) to [out=45, in=-90] (11.5-4.2,2.2);
\draw [thick] (11.2-4.2,2.5) to [out=0, in=90] (11.5-4.2,2.2);
\draw [thick] (11.2-4.2,1.7) --(9.8-4.2,.3);
\draw [thick] (9.2-4.2,1.3) --(9.8-4.2,.7);
\draw [thick] (10.2-4.2,.3) to [out=-45, in=180] (10.6-4.2,0);
\draw [thick] (11-4.2,.3) to [out=-90, in=0] (10.6-4.2,0);
\draw [thick] (11.2-4.2,2.5) --(8.8-4.2,2.5);
\draw [thick] (8.8-4.2,1.7) to [out=135, in=-90] (8.5-4.2,2.2);
\draw [thick] (8.8-4.2,2.5) to [out=180, in=90] (8.5-4.2,2.2);
\draw [thick] (11-4.2,2.7) to [out=90, in=90] (9-4.2,2.7);
\draw [thick] (11-4.2,2.3) -- (11-4.2,.3);
\draw [thick] (9-4.2,2.7) -- (9-4.2,.3);
\draw [thick] (9.8-4.2,.3) to [out=-135, in=0] (9.4-4.2,0);
\draw [thick] (9-4.2,.3) to [out=-90, in=180] (9.4-4.2,0);
\draw[fill=none] [thick](9-4.2,2.5) circle (0.2);
\draw[fill=none] [thick](11-4.2,1.5) circle (0.2);

\draw [thick] (11.2-8.3,1.7) to [out=45, in=-90] (11.5-8.3,2.2);
\draw [thick] (10.8-8.3,1.3) --(9.8-8.3,.3);
\draw [thick] (11.2-8.3,1.7) --(11.4-8.3,1.9);
\draw [thick] (9.2-8.3,1.3) --(10.2-8.3,.3);
\draw [thick] (8.8-8.3,1.7) --(8.6-8.3,1.9);

\draw [thick] (10.2-8.3,.3) to [out=-45, in=180] (10.6-8.3,0);
\draw [thick] (11-8.3,.3) to [out=-90, in=0] (10.6-8.3,0);

\draw [thick] (8.8-8.3,1.7) to [out=135, in=-90] (8.5-8.3,2.2);
\draw [thick] (11.5-8.3,2.2) to [out=90, in=90] (8.5-8.3,2.2);

\draw [thick] (11-8.3,1.9) -- (11-8.3,.3);
\draw [thick] (9-8.3,1.9) -- (9-8.3,.3);
\draw [thick] (9-8.3,1.9) to [out=90, in = 90] (11-8.3,1.9);
\draw [thick] (9.8-8.3,.3) to [out=-135, in=0] (9.4-8.3,0);
\draw [thick] (9-8.3,.3) to [out=-90, in=180] (9.4-8.3,0);
\draw[fill=none] [thick](1.7,.5) circle (0.2);

[out=0, in=90] (11.5,2.2);
\draw [thick] (10.6-12.5,1.1) --(9.8-12.5,.3);
\draw [thick] (10-12.5,2) to [out=180, in=135] (9.4-12.5,1.1);
\draw [thick] (10.6-12.5,1.1) to [out=45, in=0] (10-12.5,2);
\draw [thick] (10.2-12.5,.3) --(9.4-12.5,1.1);
\draw [thick] (10.2-12.5,.3) to [out=-45, in=180] (10.6-12.5,0);
\draw [thick] (11-12.5,.3) to [out=-90, in=0] (10.6-12.5,0);
\draw [thick] (11-12.5,1.9) -- (11-12.5,.3);
\draw [thick] (9-12.5,1.9) -- (9-12.5,.3);
\draw [thick] (9-12.5,1.9) to [out=90, in = 90] (11-12.5,1.9);
\draw [thick] (9.8-12.5,.3) to [out=-135, in=0] (9.4-12.5,0);
\draw [thick] (9-12.5,.3) to [out=-90, in=180] (9.4-12.5,0);
\draw[fill=none] [thick](-2.5,.5) circle (0.2);

\node at (-2.5,-.5){$RC \cong \Z*\Z$};
\node at (1.7,-.5){$RC \cong \Z$};
\node at (5.9,-.5){$RRC \cong \Z*\Z*\Z$};
\node at (10.1,-.5){$RRC \cong \Z*\Z$};
\end{tikzpicture}
\caption{$\Omega.2$ moves that change $RC$ or $RRC$.}
\label{twomoves}
\end{figure}
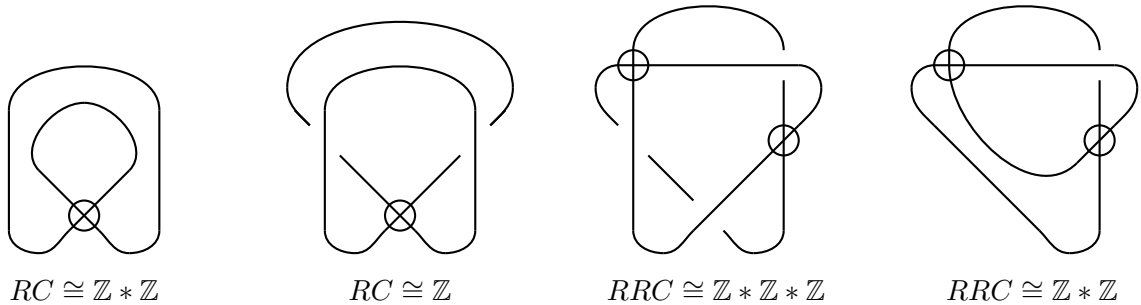

For classical link diagrams, Theorem \ref{main1} implies that $RC(D)$ is invariant under all $\Omega.2$ moves, and Theorem \ref{main2} implies that $RRC(D)$ is invariant under $\Omega.2$ moves that do not change the number $k$ of non-split pieces of $D$. In contrast, an $\Omega.2$ move on a virtual link diagram can change $RC(D)$, or change $RRC(D)$ without changing the number $k$. See Figure \ref{twomoves} for some examples. 

The following proposition shows that some $\Omega.2$ moves on virtual link diagrams do preserve $RC(D)$ and $RRC(D)$. 

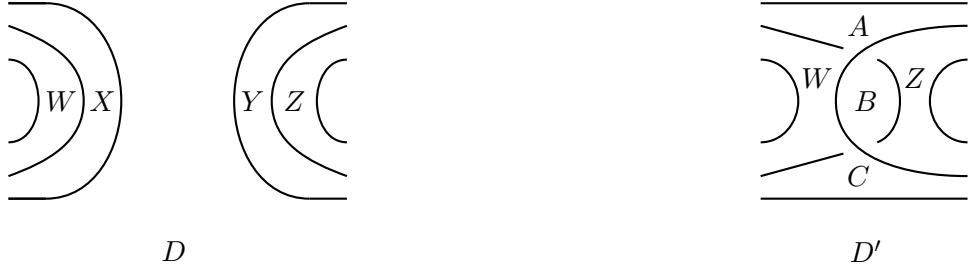
\begin{figure} 
\centering
\begin{tikzpicture} [>=angle 90]
\draw [thick] (.75,-1) to [out=-180, in=-90] (-1,0);
\draw [thick] (.75,1) to [out=-180, in=90] (-1,0);
\draw [thick] (-2,-1) -- (-.9,-0.7);
\draw [thick] (-2,1) -- (-.9,0.7);
\draw [thick] (-0.45,-0.55) to [out=20, in=-20] (-0.45,0.55);
\draw [thick] (-2,-0.55) to [out=0, in=-90] (-1.5,0);
\draw [thick] (.75,0.55) to [out=180, in=90] (.25,0);
\draw [thick] (.75,-0.55) to [out=180, in=-90] (.25,0);
\draw [thick] (-2,0.55) to [out=0, in=90] (-1.5,0);
\draw [thick] (-7.5,-0.55) to [out=180,in=-90] (-7.9,0);
\draw [thick] (-7.5,0.55) to [out=180,in=90] (-7.9,0);
\draw [thick] (-8,-1.3) to [out=180,in=-90] (-9,0);
\draw [thick] (-9,0) to [out=90,in=180] (-8,1.3);
\draw [thick] (-12,-1) to [out=20,in=-90] (-11,0);
\draw [thick] (-10.5,0) to [out=90,in=0] (-11.5,1.3);
\draw [thick] (-7.5,-1) to [out=160,in=-90] (-8.5,0);
\draw [thick] (-8.5,0) to [out=90,in=-160] (-7.5,1);
\draw [thick] (-7.5,1.3) -- (-8,1.3);
\draw [thick] (-8,-1.3) -- (-7.5,-1.3);
\draw [thick] (-12,1.3) -- (-11.5,1.3);
\draw [thick] (-12,-1.3) -- (-11.5,-1.3);
\draw [thick] (-12,-0.55) to [out=0,in=-90] (-11.6,0);
\draw [thick] (-12,0.55) to [out=0,in=90] (-11.6,0);
\draw [thick] (-11.5,-1.3) to [out=0,in=-90] (-10.5,0);
\draw [thick] (-11,0) to [out=90,in=-20] (-12,1);
\draw [thick] (-12,1.3) -- (-11.5,1.3);
\draw [thick] (-12,-1.3) -- (-11.5,-1.3);
\draw [thick] (-2,1.3) -- (.75,1.3);
\draw [thick] (-2,-1.3) -- (.75,-1.3);
\node at (-.7,1) {$A$};
\node at (-.7,-1) {$C$};
\node at (-0.6,0) {$B$};
\node at (-1.25,0.3) {$W$};
\node at (0.05,0.3) {$Z$};
\node at (-11.3,0) {$W$};
\node at (-10.75,0) {$X$};
\node at (-8.75,0) {$Y$};
\node at (-8.2,0) {$Z$};
\node at (-9.8,-2) {$D$};
\node at (-.6,-2) {$D'$};
\end{tikzpicture}
\caption{An $\Omega.2$ move changes $D$ into $D'$.}
\label{secondmove}
\end{figure}

\begin{prop}
\label{twomove}
Consider the $\Omega.2$ move pictured in Figure \ref{secondmove}. 
\begin{itemize}
\item If $\g_X=\g_Y$ in $RC(D)$, then $RC(D) \cong RC(D')$.
\item If $A$ and $C$ are distinct regions in $D'$, then $RRC(D) \cong RRC(D')$.
\end{itemize}
\end{prop}
\begin{proof} Notice that the left-hand sides of $A$ and $C$ in $D'$ come from the region $X$ of $D$, and the right-hand sides of $A$ and $C$ come from $Y$. 

Relations in $RC(D')$ include $\g_Z \g_A^{-1}=\g_W \g_B^{-1}$ and $\g_W \g_B^{-1}=\g_Z \g_C^{-1}$, implying $\g_A=\g_C$ in $RC(D')$. If $\g_X=\g_Y$ in $RC(D)$, then there is a well-defined isomorphism $f:RC(D) \to RC(D')$ with $f(\g_X)=\g_A, f(\g_W)=\g_W,f(\g_Z)=\g_Z$, and $f(\g_R)=\g_R$ for every region $R$ outside the figure.

Relations in $RRC(D')$ include $1=r_B=(g^{-1}_Bg_A)^{-1}(g^{-1}_Bg_C)$, implying $g_A=g_C$ in $RRC(D')$. If $A$ and $C$ are distinct regions in $D'$ then $X$ and $Y$ are the same region in $D$: the left-hand sides of $X$ are connected to the right-hand sides of $Y$ through $A$ on the top, and $C$ on the bottom. The relators $r_A$ and $r_C$ are $w_A(g^{-1}_Ag_B)^{-1}$ and $w_Cg^{-1}_Cg_B$, where $w_A$ and $w_B$ are words arising from the portions of the boundaries of $A$ and $C$ outside the diagram. These relators are equivalent to the equalities $g_B=g_Aw_A$ and $g_B=g_Cw^{-1}_C$. Using either equality to eliminate the generator $g_B$, we are left with the equality $g_Aw_A=g_Cw^{-1}_C$, or equivalently, $w_Aw_C=1$. As $X=Y$, $w_Aw_C$ is the same as the relator $r_X$ of $RRC(D)$. 
\end{proof}

If $D$ is a classical link diagram, then Proposition \ref{twomove} includes properties of $RC(D)$ and $RRC(D)$ that have already been mentioned. On the one hand, $X$ and $Y$ are the same region on the left-hand side of Figure \ref{secondmove}, so the hypothesis $\g_X=\g_Y$ of Proposition \ref{twomove} is satisfied automatically; thus $RC(D)$ is invariant under arbitrary $\Omega.2$ moves. On the other hand, if the two pictured sides of $D$ are not split from each other, then the two regions $A$ and $C$ of $D'$ cannot be the same, so the hypothesis $A \neq C$ of Proposition \ref{twomove} is satisfied automatically; thus $RRC(D)$ is invariant under $\Omega.2$ moves that do not change the number of split pieces of a diagram. 

\section{Abelianized core groups}
\label{abelianizations}

In this section we mention some simple properties of the abelianizations of $AC(D), RC(D)$ and $RRC(D)$. In general, if $G$ is a group we use $G_{ab}$ to denote the abelianization of $G$. 

\begin{prop}
\label{arbab}
If $A$ is an infinite, finitely generated abelian group then $A \cong AC(D)_{ab}$ for some classical link diagram $D$.
\end{prop}
\begin{proof}
For any link diagram $D$, there is an epimorphism $AC(D) \to \Z$ with $g_a \mapsto 1 \thickspace \allowbreak \forall a \in A(D)$. It follows that $AC(D)_{ab}$ is infinite.
\begin{figure} [bth]
\centering
\begin{tikzpicture}
\draw [thick] (-.4,4.6) -- (.4,5.4);
\draw [thick] (-.2,5.2) -- (-.4,5.4);
\draw [thick] (.2,4.8) -- (.4,4.6);
\draw [thick] (1.1,4.6) -- (1.9,5.4);
\draw [thick] (1.3,5.2) -- (1.1,5.4);
\draw [thick] (1.7,4.8) -- (1.9,4.6);
\draw [thick] (.4,4.6) to [out=-45, in=-135] (1.1,4.6);
\draw [thick] (.4,5.4) to [out=45, in=135] (1.1,5.4);
\draw [thick] (2.6,4.6) -- (3.4,5.4);
\draw [thick] (2.8,5.2) -- (2.6,5.4);
\draw [thick] (3.2,4.8) -- (3.4,4.6);
\draw [thick] (1.9,4.6) to [out=-45, in=-135] (2.6,4.6);
\draw [thick] (1.9,5.4) to [out=45, in=135] (2.6,5.4);
\draw[fill=black] [thick](4,5) circle (.03);
\draw[fill=black] [thick](4.2,5) circle (.03);
\draw[fill=black] [thick](3.8,5) circle (.03);
\draw [thick] (2+2.6,4.6) -- (2+3.4,5.4);
\draw [thick] (2+2.8,5.2) -- (2+2.6,5.4);
\draw [thick] (2+3.2,4.8) -- (2+3.4,4.6);
\draw [thick] (5+.4,4.6) to [out=-45, in=-135] (5+1.1,4.6);
\draw [thick] (5+.4,5.4) to [out=45, in=135] (5+1.1,5.4);
\draw [thick] (3.5+2.6,4.6) -- (3.5+3.4,5.4);
\draw [thick] (3.5+2.8,5.2) -- (3.5+2.6,5.4);
\draw [thick] (3.5+3.2,4.8) -- (3.5+3.4,4.6);
\draw [thick] (6.9,5.4) to [out=45, in = 135] (-.4,5.4);
\draw [thick] (6.9,4.6) to [out=-45, in = -135] (-.4,4.6);
\end{tikzpicture}
\caption{This diagram of a $(2,m)$ torus link has $m$ crossings.}

\label{toruslink}
\end{figure}
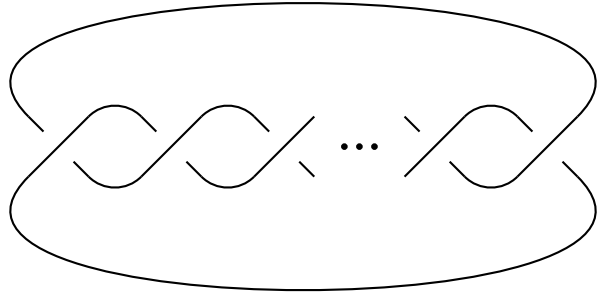

The reader can easily verify that if $D$ is a classical link diagram, $m \in \mathbb N$ and $D'$ is the connected sum of $D$ with the $(2,m)$ torus link diagram pictured in Figure \ref{toruslink}, then $AC(D')_{ab} \cong AC(D)_{ab} \oplus \Z_m$. If $A$ is an infinite, finitely generated abelian group, then 
\[
A \cong \Z^p  \oplus \Z_{m_1} \oplus \cdots \oplus \Z_{m_q}
\]
for some choice of integers $p \geq 1$, $q \geq 0$ and $m_1, \dots, m_q \geq 2$. We can construct a classical link diagram $D$ with $AC(D)_{ab} \cong A$ by starting with $p$ disjoint circles, and then taking connected sums with torus links of types $(2,m_1) , \dots, (2,m_q)$.\end{proof}

 \begin{prop}
Let $D$ be a diagram of a $\mu$-component link. Then the 2-rank of $AC(D)_{ab}$ is $\mu$.
 \end{prop}
 \begin{proof}
The abelian version of the crossing relation $1=g_{a_1(c)}g^{-1}_{a_2(c)}g_{a_1(c)}g^{-1}_{a_3(c)}$ of Definition \ref{arccore} is $1=g^2_{a_1(c)}g^{-1}_{a_2(c)}g^{-1}_{a_3(c)}$. If we reduce exponents modulo 2 we obtain $1=g^{-1}_{a_2(c)}g^{-1}_{a_3(c)}$, which is equivalent to $g_{a_2(c)}=g_{a_3(c)}$. That is, the quotient group $AC(D)_{ab}/(2 \cdot AC(D)_{ab})$ has a generator for each arc of $D$, with relations that imply that at each crossing, the generators corresponding to the two underpassing arcs are equal. It follows that for each component of the link that appears in the diagram, all the generators of $AC(D)_{ab}/(2 \cdot AC(D)_{ab})$ corresponding to arcs of that component are equal.  \end{proof}

Notice that Theorems \ref{main1} and \ref{main2} imply that for a classical link diagram $D$, if $T$ is the torsion subgroup of $AC(D)_{ab}$ then $T$ and $T \oplus T$ are isomorphic to the torsion subgroups of $RC(D)_{ab}$ and $RRC(D)_{ab}$, respectively. The examples in Figure \ref{virttref} imply that in contrast, for virtuals the three torsion groups are independent.

\section{Another description of $AC(D)$}\label{another}

In this section we discuss the fact that in addition to  Definition \ref{arccore}, the arc core group $AC(D)$ has another, somewhat more complicated description. Our discussion is modified from J. Przytycki's account \cite{P}.

Let $D$ be an oriented virtual link diagram, $D^+$ the diagram obtained by adjoining an unknotted circle with no crossings. Let $F_{A(D^+)}$ be the free group on the set $\{x_a \mid a \in A(D^+) \}$, and let $F_{A(D^+)}^{(2)}$ be the subset of $F_{A(D^+)}$ consisting of elements $\prod x_{a_i}^{n_i}$ such that $\sum n_i$ is even. The sum $\sum n_i$ is not changed by conjugation, so $F_{A(D^+)}^{(2)}$ is a normal subgroup of $F_{A(D^+)}$. Let $N$ be the normal subgroup of $F_{A(D^+)}^{(2)}$ generated by $\{ x^2_a \mid a \in A(D^+) \}$.

For each crossing $c \in C(D)$, let $a_1(c)$ be the overpassing arc at $c$ and $a_2(c),a_3(c)$ the underpassing arcs at $c$, labeled in such a way that an observer standing on $c$ facing forward along $a_1(c)$ has $a_2(c)$ on the right. The Wirtinger relator $w(c)$ is then $w(c)=x_{a_1(c)}x_{a_2(c)}x^{-1}_{a_1(c)}x^{-1}_{a_3(c)} 
$. Let $\widehat W$ be the normal subgroup of $F_{A(D^+)}^{(2)}/N$ generated by $\{w(c)N \mid c \in C(D)\}$. 

\begin{lemma}
\label{normalclosure}
$N$ is a normal subgroup of $F_{A(D^+)}$, and $\widehat W$ is a normal subgroup of $F_{A(D^+)}/N$.
\end{lemma}
\begin{proof}
To show that $N$ is normal in $F_{A(D^+)}$, it suffices to show that for every choice of $a,b \in A(D^+)$, $x_b x_a^2 x^{-1}_b$ and $x^{-1}_b x^2_a x_b$ are both elements of $N$. As $x_b x^2_a x^{-1}_b = (x_bx_a) x^2_a (x_b x_a)^{-1}$ and $
x^{-1}_b x^2_a x_b = (x_ax_b)^{-1} x^2_a (x_a x_b)$, $x_b x^2_a x^{-1}_b$ and $x^{-1}_b x^2_a x_b$ are both equal to conjugates of $x^2_a$ by elements of $F_{A(D^+)}^{(2)}$. As $N$ is normal in $F_{A(D^+)}^{(2)}$, it follows that $x_b x^2_a x^{-1}_b$ and $x^{-1}_b x^2_a x_b$ are both elements of $N$.

Similarly, to show that $\widehat W$ is normal in $F_{A(D^+)}/N$ we must show that for every choice of $b \in A(D^+)$ and $c \in C(D)$, $x_b w(c) x^{-1}_b N$ and $x^{-1}_b w(c) x_b N$ are both elements of $\widehat W$. We begin to verify this by observing that
\[
w(c)^{-1} N = x_{a_3(c)} x_{a_1(c)}x^{-1}_{a_2(c)}x^{-1}_{a_1(c)}N=x^2_{a_3(c)} x^{-1}_{a_3(c)} x_{a_1(c)}  x^{-2}_{a_2(c)} x_{a_2(c)}x^{-1}_{a_1(c)}N
\]
\[
= x^2_{a_3(c)}N \cdot x^{-1}_{a_3(c)} x_{a_1(c)}N \cdot x^{-2}_{a_2(c)}N \cdot x_{a_2(c)}x^{-1}_{a_1(c)}N 
\]
\[
= 1N \cdot x^{-1}_{a_3(c)} x_{a_1(c)}N \cdot 1N \cdot x_{a_2(c)}x^{-1}_{a_1(c)}N = x^{-1}_{a_3(c)} x_{a_1(c)} x_{a_2(c)}x^{-1}_{a_1(c)}N = x^{-1}_{a_3(c)} w(c) x_{a_3(c)}N.
\]
Taking inverses, we see that $w(c)N = x^{-1}_{a_3(c)} w(c)^{-1} x_{a_3(c)}N$. 

The fact that $N$ is normal in $F_{A(D^+)}$ implies that $F_{A(D^+)}^{(2)}/N$ is simply a subgroup of $F_{A(D^+)}/N$. Therefore we can perform calculations involving elements of $F_{A(D^+)}^{(2)}/N$ in $F_{A(D^+)}/N$; the results will be the same. It follows that $w(c)N = x^{-1}_{a_3(c)} w(c)^{-1} x_{a_3(c)}N$ implies
\[
x_b w(c) x^{-1}_b N = x_b x^{-1}_{a_3(c)} w(c)^{-1} x_{a_3(c)} x^{-1}_b N = (x_b x^{-1}_{a_3(c)}N) \cdot w(c)^{-1}N \cdot (x_b x^{-1}_{a_3(c)}N)^{-1}
\]
and similarly,
\[
x^{-1}_b w(c) x_b N = x^{-1}_b x^{-1}_{a_3(c)} w(c)^{-1} x_{a_3(c)} x_b N= (x^{-1}_b x^{-1}_{a_3(c)}S) \cdot  w(c)^{-1}N \cdot (x^{-1}_b x^{-1}_{a_3(c)} N)^{-1}.
\]
Therefore each of $x_b w(c) x^{-1}_b N,x^{-1}_b w(c) x_b N$ is equal to a conjugate of $w(c)^{-1}N$ by an element of $F_{A(D^+)}^{(2)}/N$. As $w(c)^{-1}N \in \widehat W$ and $\widehat W$ is normal in $F_{A(D^+)}^{(2)}/N$, it follows that $x_b w(c) x^{-1}_b N$ and $x^{-1}_b w(c) x_b N$ are both elements of $\widehat W$. \end{proof}

As $F_{A(D^+)}$ is a free group, its subgroup $F^{(2)}_{A(D^+)}$ is also free. If $a^+ \in A(D^+)$ is the arc of $D^+$ representing the extra component added to $D$, then $F^{(2)}_{A(D^+)}$ is freely generated by the set $\{x_{a^+}^2\} \cup \{ x^2_a,x_ax_{a^+} \mid a \in A(D) \}$. It follows that the quotient group $F^{(2)}_{A(D^+)}/N$ is freely generated by $\{ x_ax_{a^+}N \mid a \in A(D) \}$. Therefore if $F_{A(D)}$ is the free group on $\{g_a \mid a \in A(D)\}$, there is an isomorphism $f:F_{A(D)} \to F^{(2)}_{A(D^+)}/N$  with $f(g_a)=x_ax_{a^+}N \thickspace \allowbreak \forall a \in A(D)$.
\begin{theorem}
\label{secondpres}
The normal subgroup of $F_{A(D)}$ generated by the relators $r_c$ of Definition \ref{arccore} is $f^{-1}(\widehat W)$.
\end{theorem}
\begin{proof}
The theorem follows immediately from the fact that 
\[
f(r_c)=f(g_{a_1(c)}g^{-1}_{a_2(c)}g_{a_1(c)}g^{-1}_{a_3(c)})=x_{a_1(c)}x^{-1}_{a_2(c)}x_{a_1(c)}x^{-1}_{a_3(c)}N
\]
\[
=x_{a_1(c)}x^{-1}_{a_2(c)}N \cdot 1N \cdot x_{a_1(c)}x^{-1}_{a_3(c)}N =x_{a_1(c)}x^{-1}_{a_2(c)}N \cdot x^{2}_{a_2(c)}x^{-2}_{a_1(c)}N \cdot x_{a_1(c)}x^{-1}_{a_3(c)}N
\]
\[
=x_{a_1(c)}x_{a_2(c)}x^{-1}_{a_1(c)}x^{-1}_{a_3(c)}N =w(c)N.
\]
\end{proof}

\begin{cor}
\label{coreiso}
The map $f$ defines an isomorphism $AC(D) \cong F^{(2)}_{A(D^+)}/\langle NW \rangle$, where $\langle NW \rangle$ is the normal subgroup of $F_{A(D^+)}^{(2)}$ generated by $\{ x^2_a \mid a \in A(D^+) \} \cup \{ w(c) \mid c \in C(D)\}$. 
\end{cor}
\section{Core group functor}  It is well known that the arc core group of a classical link is related to its Wirtinger group. One might anticipate a similar connection between the region core group and the Dehn group. We will explain both relationships, extending the first to the virtual case.

We recall that the \emph{Wirtinger group} $\pi_{wirt}(D)$ of an oriented virtual link diagram $D$ has a presentation with generators corresponding to the arcs of $D$. Relators correspond to each classical crossing $c$ as follows: if $a_1(c)$ is the overpassing arc at $c$ and $a_2(c),a_3(c)$ the underpassing arcs at $c$, labeled in such a way that an observer standing on $c$ facing forward along $a_1(c)$ has $a_2(c)$ on the right, then the Wirtinger relator corresponding to $c$ is $x_{a_1(c)}x_{a_2(c)}x^{-1}_{a_1(c)}x^{-1}_{a_3(c)}$. The choice of orientations of the components of $D$ does not affect $\pi_{wirt}(D)$. It is well known that $\pi_{wirt}(D)$ is invariant under virtual Reidemeister moves, and hence it is an invariant $\pi_{wirt}(L)$ of the virtual link $L$ described by $D$. 

The \emph{Dehn group} $\pi_{dehn}(D)$ is another group that can be associated to $D$, although no orientation is necessary. It has a presentation with generators corresponding to the regions of $D$, with one generator set equal to the identity, and relators at each crossing of the form $\x_V^{}\x_W^{-1}\x_X^{} \x_Y^{-1}$, as in Figure \ref{regionfig}. 

\begin{remark} (1.) It is well known that when $L$ is a classical link,  the Wirtinger and Dehn groups are both isomorphic to the fundamental group of the link complement $\R^3 \setminus L$. Although the Wirtinger group is also a virtual link invariant, the Dehn group is not. 

(2.) If $L$ is regarded as a link in a thickened surface, then both $\pi_{wirt}(D)$ and $\pi_{dehn}(D)$ are isotopy invariants. (The proof as in the classical case, using Reidemeister moves, extends.) However, the groups are in general no longer isomorphic. In fact, $\pi_{dehn}(D)$ is a quotient of $\pi_{wirt}(D)$ by a normal subgroup that is easy to describe. We review this below. 

(3.) It is a result of \cite{BH} that any finitely presented group arises as the Dehn group of some virtual link. \end{remark} 

Let $G$ be a finitely presented group, $S \subset G$ a finite subset and $N$ the normal closure of $\{s^2 \vert s \in S\}$. 
Assume that 
$\rho:  G \to \Z_2$ is a homomorphism that maps each $s\in S$ nontrivially. 
Define $G^+$ to be the free product $G*\<y\>$  and $N^+$ the normal closure of  $\{s^2 \vert s \in S\} \cup \{y^2\}$.
Extend $\rho$ to a homomorphism $\rho^+: G^+ \to \Z_2$ by mapping $y$ nontrivially.  

\begin{definition} \label{defcore}  $\C(G, S, \rho)$ is the group ${\rm Ker} \rho^+/ N^+$. \end{definition}

\begin{remark} One checks that $\C$ is a covariant functor from the category {\bf Core} to the category of groups.  Objects of ${\bf Core}$ are 3-tuples $(G, S, \rho)$. Morphisms $F: (G_1, S_1, \rho_1) \to 
(G_2, S_2, \rho_2)$ are homomorphisms $F: G_1 \to G_2$  such that $F(N_1) \subset N_2$ and 
$\rho_2 \circ F = \rho_1$.  \end{remark}

\begin{definition} \label{standard} Let  $D$ be an oriented virtual link diagram.  The \emph{standard (Wirtinger) 3-tuple} of $D$ is $(\pi_{wirt}(D), S, \rho)$, where $\pi_{wirt}(D)$ is the Wirtinger group of the diagram, $S$ is the set of Wirtinger generators $x_i$, and $\rho: \pi_{wirt}(D) \to \Z_2$ is the unique homomorphism that maps each $x_i$ nontrivially. \end{definition}

Consider a standard 3-tuple $(\pi_{wirt}(D), S, \rho)$ of a diagram $D$ of a virtual link $L$. Define $L^+$ to be the virtual link comprising $L$ and a distant unknotted component with meridanal generator $y$. Any virtual link diagram $D^+$ of $L^+$ can be converted to another virtual link diagram of $L^+$ by a finite sequence of generalized Reidemeister moves. The groups presented by the Wirtinger presentations of equivalent virtual link diagrams are isomorphic by an isomorphism that preserves conjugacy classes of Wirtinger generators. Covariance of the core group functor $\C$ ensures that $\C(\pi_{wirt}(D), S, \rho)$ is independent of the diagram of $L$. It is also independent of the diagram's orientation. We will abbreviate the group $\C(\pi_{wirt}(D), S, \rho)$ by  $\C(\pi_{wirt}(D))$.

\begin{theorem}\label{catcore} If $L$ is a virtual link, then $\C(\pi_{wirt}(D)) \cong AC(D)$. \end{theorem}

\begin{proof} A strictly algebraic proof of this theorem is given in Section \ref{another}. Other methods such as the Reidemeister-Schreier method also can be used. Here we give a topological argument, using covering spaces.

Assume that $D$ is an oriented virtual link diagram for $L$ with Wirtinger generators $x_a$, for $a \in A(D)$, and
Wirtinger relators of the form $x_{a_1(c)}x_{a_2(c)}x^{-1}_{a_1(c)}x^{-1}_{a_3(c)}$ at each crossing $c$. Consider the standard 2-complex associated to the presentation. It has a single 0-cell $b$, oriented 1-cells $\s_a$ corresponding to the Wirtinger generators $x_a$, and 2-cells
$\t_c$ with attaching maps described by the relators. Let $X$ be the complex with an additional 1-cell $\s_y$ to realize the Wirtinger group of $D^+$. 

Let $\widetilde X$ be the 2-fold covering space with fundamental group 
${\rm Ker}\ \rho^+$. It has a lifted complex structure consisting of 0-cells $\ti b, \ti b'$;  oriented 1-cells $\ti \s_a, \ti \s_y$ (resp. $\ti \s'_a, \ti \s'_y$) from $\ti b$ to $\ti b'$ (resp. from $\ti b'$ to $\ti b$);  2-cells $\ti \t_c$ (resp. $\ti \t'_c$) that are lifts of the $\t_c$ with boundaries based at $\ti b$ (resp. $\ti b'$).
The fundamental group $\pi_1(\widetilde X, \ti b)$ is generated by the homotopy classes of closed paths based at $\ti b$, the lifts of closed paths in $X$ based at $b$ of even word-length. 

For each 1-cell $\s$ in $X$ attach a 2-cell to $\widetilde X$ along the closed path $\ti \s \ti \s'$. Denote the new cell complex by $\widehat X$. Then the normal subgroup in $\pi_1(\widetilde X, \ti b)$ generated by the classes of these closed paths is the subgroup $N^+$ of Definition \ref{defcore}. Hence $\pi_1(\widehat X, \ti b) \cong \C(\pi_{wirt}(D))$. 

We complete the proof by showing that $\pi_1(\widehat X, \ti b)$ is isomorphic to $AC(D)$. To do this, contract the closed path $\ti \s_y \ti \s'_y$ together with its bounding 2-cell, identifying $\ti b'$ with $\ti b$ and 
converting all remaining 1-cells in $\widehat X$ to oriented loops based at $\ti b$. We abuse notation by denoting the converted 1-cells (now loops) by the same symbols as before.  Contracting $\ti \s_y \ti \s'_y$ and the bounding 2-cell preserves the homotopy type and hence the fundamental group of the complex $\widehat X$. 

We denote the homotopy class of a based oriented closed path with square brackets $[\  \cdot \ ]$. Note that  $[\ti \s'_a] =
[\ti\s_a]^{-1}$. The boundaries of the remaining 2-cells $\ti \t_c, \ti \t_c'$ can be read as 
$\ti\s_{a_1(c)}\ti\s_{a_2(c)}^{-1} \ti\s_{a_1(c)}\ti\s_{a_3(c)}^{-1}$ and $\ti\s_{a_1(c)}^{-1}\ti\s_{a_2(c)} \ti\s_{a_1(c)}^{-1}\ti\s_{a_3(c)}$, respectively. (Note the alternating exponents.)
The group elements  $[\ti\s_{a_1(c)}\ti\s_{a_2(c)}^{-1} \ti\s_{a_1(c)}\ti\s_{a_3(c)}^{-1}]$ and $[\ti\s_{a_1(c)}^{-1}\ti\s_{a_2(c)} \ti\s_{a_1(c)}^{-1}\ti\s_{a_3(c)}]^{-1}$ are conjugates, so we can ignore the 2-cells $\ti\t_j'$ without affecting the fundamental group. We obtain the presentation of $AC(D)$ in Definition \ref{arccore}.   Hence $\C(\pi_{wirt}(D)) \cong \pi_1(\widehat X, \ti b)$ is isomorphic to $AC(D)$. \end{proof}

The above proof suggests an easy procedure for finding a presentation of $\C(G, S, \rho)$ when $S$ generates $G$. We describe it next. 

\begin{definition} If $w = g_{i_1}^{\e_1} g_{i_2}^{\e_2} \cdots  g_{i_k}^{\e_k}$ is a word in the free group 
generated by $g_1, \ldots, g_n$ and $\e_i =\pm 1$, then $\ti w=
g_{i_1}g_{i_2}^{-1}g_{i_3} \cdots g_{i_k}^{(-1)^{k-1}}$ and $(\ti w)' =   g_{i_1}^{-1}g_{i_2} g_{i_3}^{-1} \cdots  g_{i_k}^{(-1)^k}.$\end{definition}

From the proof of Theorem \ref{catcore} we have: 
\begin{prop}\label{compute} Assume that $(G, S, \rho)$ is a triple in the category {\bf Core} such that $G$ has a presentation of the form 
$\<S ; r_1, \ldots, r_m\>$. Then $\C(G, S, \rho)\cong \<S; \ti r_1, (\ti r_1)', \ldots,  \ti r_m, (\ti r_m)' \>$. \end{prop}

\begin{remark} (1.) A generator $s \in S$ of  $G$ should not be confused with the generator $s \in \C(G, S, \rho)$. However, the two have a simple relationship. Consider the cover $\widehat X $ before the 1-cell $\s'_y$ is contracted. A generator $s \in G$ is represented by a loop in $X$ that lifts to directed 1-cell $\s_s$ in $\widehat X$. Together with $\s'_y$, it determines a directed loop in $\widehat X$ representing $s \in \C(G, S, \rho)$. 

(2.) It is often the case that the relator $(\ti r_i)'$ is a consequence of $\ti r_i$ and hence can be omitted from the presentation in Proposition \ref{compute}.  \end{remark}

\begin{cor} A group is isomorphic to the core group of a virtual knot if and only if it has a presentation of the form 
$$\< x_1, \ldots, x_n ; \T{u_i x_{i+1}} = \T{x_i u_i}\ (i \in \Z_n) \>,$$
where $u_i$ are arbitrary elements of the free group on $x_1, \ldots, x_n$. \end{cor}

\begin{proof} Theorem 2.2 of  \cite{SW00} states that a group is the Wirtinger group of a virtual knot if and only if it has a presentation of the form $\< x_1, \ldots, x_n ; u_i x_{i+1} = x_i u_i\ (i \in \Z_n) \>$, where $u_i$ are elements of the free group on $x_1, \ldots, x_n$. 
Consider the standard core 3-tuple of such a group. Its core group has a presentation 
$$\< x_1, \ldots, x_n ; \T{u_i x_{i+1}} = \T{x_i u_i}, (\T{u_i x_{i+1}})' = (\T{x_i u_i})'\ (i \in \Z_n) \>.$$
It is easily seen that the relations $(\T{u_i x_{i+1}})' = (\T{x_i u_i})'$ are redundant and so can be omitted. 
\end{proof}

The Wirtinger group of any classical 2-bridge knot has a simplified presentation $\< x_1, x_2 ; r\>$, where $r$ has the form $ux_2 = x_1u$, for some word $u$ in $x_1, x_2$. 
Then $\ti r$ can be seen to have the form $(x_1x_2^{-1})^d$, for some positive integer $d$. The following well-known result is immediate. 

\begin{cor} \label{2bridge} The core group of any classical 2-bridge knot is isomorphic to $\Z * \Z_d$ for some positive integer $d$. \end{cor}

The integer $d$ in Corollary \ref{2bridge} is the {\sl determinant} of the knot, the order of the fundamental group of its 2-fold branched cyclic cover.

\begin{remark}\label{defvar}  (1.) The proof of Theorem \ref{catcore} suggests alternative ways to compute $\C(G, S, \rho)$: In Definition \ref{defcore} we extend $G$ to $G^+$ in order to achieve the desired symmetry in the presentation of the core group. However, we can obtain $\C(G, S, \rho)$ without extending $G$. In the construction of the proof of Theorem \ref{catcore}, we can omit the 1-cell $y$ in the construction of the complex $X$. Then instead of contracting the 1-cells $\s_y, \s'_y$ and bounding 2-cell, we can choose any $s \in S$ instead of $y$, and contract $\s_s \cup \s'_s$ and the bounding 2-cell. The group $\C(G, S, \rho)$ is isomorphic to $\Z * {\rm Ker} \rho /N$.

(2.) Since the homomorphism $\rho$ vanishes on the subgroup $N$,  it induces a well-defined homomorphism $\bar \rho: G/N \to \Z_2$. Using the previous remark, we see that  $\C(G, S, \rho) \cong \Z* {\rm Ker} \bar\rho$. \end{remark}

\begin{prop} \label{braidcore} Let ${\cal B}_n$ be the $n$-string braid group, $n\ge 2$, with standard generators $b_1, \dots, b_{n-1}$. Let $S=\{b_1, \ldots, b_{n-1}\}$ and $\rho: {\cal B}_n \to \Z_2$ the unique homomorphism mapping each $b_i$ nontrivially. Then $\C({\cal B}_n, S, \rho)$ is isomorphic to $\Z*A_n$, where $A_n$ is the alternating group on $n$ symbols. \end{prop}

\begin{proof} We recall that ${\cal B}_n$ has presentation 
$$\< b_1, \ldots, b_{n-1} ; b_i b_{i+1} b_i = b_{i+1} b_i b_{i+1}\ (\forall i), b_i b_j = b_j b_i \ (\forall i, j, |i-j|>1) \>.$$
The assignment of each $b_i$ to the transposition $s_i = (i\ i+1)$ induces a homomorphism from ${\cal B}_n$ onto the symmetric group $S_n$. 
Its kernel $P{\cal B}_n$, the {\it pure braid group}, was shown by Artin to be the normal subgroup of ${\cal B}_n$ generated by 
$b_1^2, \ldots, b_{n-1}^2$. (See, for example, Chapter 9 of \cite{FM}.)  The kernel of the induced homomorphism $\bar \rho: {\cal B}_n/P{\cal B}_n \to \Z_2$ is isomorphic to $A_n$, the subgroup of $S_n$ consisting of words of even length in the generators $b_i$.  By Remark \ref{defvar}, $\C({\cal B}_n, S, \rho)$ is isomorphic to $\Z*A_n$. 
 
\end{proof} 

It is well known that the alternating group has a presentation 
$$A_n = \< a_1, \ldots, a_{n-2}; a_i^3, (a_i a_j)^2 \ (\forall i, j, 1 \le i \ne j \le n-2)\>$$
(see p. 64 of \cite{Bo}). The procedure from the proof of Theorem \ref{catcore} for finding a presentation of $\C(G, S, \rho)$, contracting $\s_{b_{n-1}}, \s'_{b_{n-1}}$ and bounding 2-cell, yields a similar presentation for $A_n$. 

\begin{cor} \label{alt} The alternating group $A_n$ has presentation 
$$\< x_1, \ldots, x_{n-1} ; (x_i x^{-1}_{i+1})^3 \ (\forall i), (x_i x^{-1}_j)^2 \ (\forall i,j, |i-j|>1), x_{n-1} \>.$$
\end{cor}

\begin{example} By  Corollary \ref{alt}, we obtain
$$A_5 \cong  \< x_1, x_2, x_3; (x_1x_2^{-1})^3, (x_2 x_3^{-1})^3, x_3^3, (x_1x_3^{-1})^2, x_1^2, x_2^2 \>.$$
Using Tietze transformations, we introduce new generators $A, B, C$ and defining relations $A = x_1x_2^{-1}, B= x_2 x_3^{-1}, C = x_3$. Then
$$A_5 \cong \< A, B, C;  A^3=B^3=C^3 = (AB)^2= (BC)^2=(ABC)^2=1 \>.$$
(Here $A = s_1s_2, B= s_2s_3$ and $C = s_3s_4$.) 
The presentation of $A_5$, originally discovered by J.A. Todd in 1931, appears on p. 125 of \cite{CM}: 
$$A_5 \cong \< V_1, V_2, V_3;  V_1^3=V_2^3=V_3^3 = (V_1V_2)^2= (V_2V_3)^2=(V_3V_1)^2=1 \>.$$
\end{example}


\section{Coxeter groups}

Let $M = (m_{ij}), 1 \le i, j \le n$, be a symmetric $n \times n$ matrix with entries in ${\mathbb N} \cup \{\infty\}$, such that $m_{ii}=1$ for all $i$ and 
$m_{ij}>1$ if $i \ne j$. The \emph{Coxeter group of type $M$} (see \cite{T}) is the group defined by the presentation 
\begin{equation}\label{coxeter} W= \<s_1, \ldots, s_n ; (s_i s_j)^{m_{ij}}\ (\forall i, m_{ij}< \infty)\>. \end{equation}
Since diagonal entries of $M$ are equal to $1$, each generator $s_i$ satisfies $s_i^2 =1$. Consequently, $m_{ij}=2$ implies that $s_i$ and $s_j$ commute.

The information of the matrix $M$ can be encoded by a  graph with $n$ vertices labeled $s_1, \ldots, s_n$. There is an edge connecting $s_i$ and $s_j$ and labeled by $m_{ij}$ whenever $m_{ij}>2$. By convention, edges labeled  2 are omitted, while any unlabeled edge is assumed to be labeled 3.  

Any Coxeter group $W$ has a natural associated triple $(W, S, \rho)$, where $S=\{s_1, \ldots, s_n\}$ and $\rho: W \to Z/2$ is the unique homomorphism mapping each $s_i$ nontrivially. We will call $(W, S, \rho)$ the \emph{standard triple} of $W$. We denote by $\C(W)$
the group determined by the core functor with the standard triple of $W$.

\begin{prop}\label{coxetercompute} If $W$ is the Coxeter group of type $M$, then $\C(W)$ is isomorphic to the free product of $\Z$ and 
the kernel of $\rho$.  \end{prop}

\begin{proof} Since each square $s_i^2$ is trivial in $W$, the subgroup $S$ is trivial. By Remark 2 of  \ref{defvar},  $\C(W, S, \rho)$ is isomorphic to $\Z* {\rm Ker}\rho$. \end{proof}

\begin{example}  The symmetric group $S_n$ generated by transpositions $s_i = (i\ i+1), 1\le i <n$, is the Coxeter group corresponding to the graph that consists of $n-1$ vertices joined by $n-2$ unlabeled edges. Consider the standard triple $(S_n, S, \rho)$. The kernel of $\rho$ is the alternating group $A_n$. Hence $\C(S_n) \cong \Z*A_n$. (Compare with Proposition \ref{braidcore}.)
\end{example}

\begin{example} The group of the $(2,m)$-torus link has a presentation $\pi_{wirt} \cong \<s_1, s_2 ; s_1s_2s_1\cdots =s_2s_1s_2\cdots\>$, where each side of the relation has length $m$. Both $s_1$ and $s_2$ are Wirtinger generators, and in fact the presentation is obtained from a Wirtinger presentation by applying Tietze transformations. Consider the standard triple $(\pi_{wirt}, S, \rho)$. Then the quotient group $\pi_{wirt}/N$ has presentation $\<s_1, s_2 ; s_1^2, s_2^2, (s_1s_2)^m\>$. This is a Coxeter group corresponding to the graph consisting of two vertices joined by a single edge labeled $m$. 

Proposition \ref{compute} can be used to see that $\C(\pi_{wirt})$ has presentation $\<s_1, s_2; (s_1s_2^{-1})^m\>$. Let $U=s_1s_2^{-1}$. Then Tietze transformations convert the presentation to $\<s_1, U; U^m\> \cong \Z* \Z_m$.

\end{example}

\begin{example} Consider the Coxeter group $W\cong \<s_1, s_2, s_3; s_1^2, s_2^2, s_3^2, (s_1s_2)^3, (s_2s_3)^3, (s_3s_1)^3\>$ corresponding to the graph consisting of a triangle with unlabeled edges. Consider the standard triple $(W, S, \rho)$.  The method of the proof of Theorem \ref{catcore} can be used to see that $\C(W)$ has presentation $\<s_1, s_2, s_3; (s_1s_2^{-1})^3, (s_2 s_3^{-1})^3, (s_3s_1^{-1})^3 \>$. Let $U= s_1s_2^{-1}$ and $V= s_2s_3^{-1}$. Then Tietze transformations convert the presentation to $\<s_1, U, V; U^3, V^3, (UV)^3\> \cong \Z * \< U, V; U^3, V^3, (UV)^3\>$. 

The group $W$ is the plane crystallographic (``wallpaper") group {\bf p31m} while $\C(W)$ is the free product of $\Z$ and the index-2 subgroup {\bf p3} of {\bf p31m} (see pp. 48--49 of \cite{CM}). 
\end{example}


\section{Dehn core groups from Dehn groups}
\label{Dehnc}
We  return our attention to the Dehn group of a virtual link diagram $D$ and its relationship to the region core group $RC(D)$.

Recall that $D$ corresponds to an {\sl abstract link diagram}, a classical link diagram on a surface $\Si$
for which the underlying 4-valent graph is a deformation retract. By attaching 2-disks along the boundary of the surface we obtain a link diagram in a closed surface ${\cal S}$. The regions of the abstract link diagram are preserved (although they are enlarged by the added disks), and so we can find the presentation of $\pi_{dehn}(D)$ from the diagram $D \subset {\cal S} $. It is straightforward to see that $\pi_{dehn}(D)$ is unchanged by Reidemeister moves in ${\cal S}$ as well as surface automorphisms. If we thicken ${\cal S}$, regarding it as ${\cal S} \times \{0\} \subset {\cal S}\times [-1, 1]$, then $D$ describes a link $L \subset {\cal S}\times [-1,1]$, and $\pi_{dehn}(D)$ is an invariant of $L$ up to ambient isotopy and automorphisms of ${\cal S} \times [-1, 1]$. Conversely, it is straightforward to see that any link $L$ in a thickened surface ${\cal S} \times [-1,1]$ arises in the above manner from a diagram $D \subset {\cal S}$. 
Henceforth surfaces ${\cal S}$ are assumed to be closed, connected and orientable.

A link diagram $D \subset {\cal S}$ is checkerboard colorable if and only if $L$ represents the trivial element in the mod-2 homology of ${\cal S}\times [-1,1]$. 
(This can be proved by observing that the intersection of $D$ with each boundary arc of a fundamental region $\De$ of ${\cal S}$ has an even number of points if and only if $D$ is checkerboard colorable.) Recall that $\pi_{dehn}(D)$ has Dehn 
generators $\x_S$ corresponding to shaded regions $S$ together with Dehn generators $\x_U$ corresponding to unshaded regions $U$, with a single shaded generator $\x_{S_0}$ set equal to the identity. Lemma \ref{grouplem} and Corollary \ref{interchange} tell us that it does not matter which shaded region is chosen as $S_0$. (One can equally well chose the region to be unshaded.) When $D$ is checkerboard colorable,  there is a unique homomorphism $\rho: \pi_{dehn}(D) \to \Z_2$  mapping each $\x_S$ to 0 and each $\x_U$ to 1.

Each arc of $D$ determines elements of the form $\x_U^{}  \x_S^{-1}$, with $U$ on one side of the arc and $S$ on the other side (at the same point). We will refer to them as {\sl arc elements} of $\pi_{dehn}(D)$. We define a {\sl standard Dehn 3-tuple} of $D$ to be $(\pi_{dehn}(D), S, \rho)$, where $S$ is the set of arc elements. Let $\rho: \pi_{dehn}(D) \to \Z_2$ be the homomorphism described above, sending each $\x_S$ to 0 and $\x_U$ to 1. Note that the 3-tuple depends on the choice of a generator that we set equal to the identity in the definition of the Dehn group, a choice that determines the homomorphism $\rho$. When the choice of generator is clear, we refer to $(\pi_{dehn}(D), S, \rho)$ as {\sl the} standard Dehn 3-tuple and shorten the notation $\C(\pi_{dehn}(D), S, \rho)$ to $\C(\pi_{dehn}(D))$. 
 
Assume that $\x_{S_0}$ is the generator set equal to the identity element in the definition of the Dehn group. By Proposition \ref{freeproduct}, the element $\x_{S_0}$ generates an infinite cyclic free factor of $RC(D)$. The remaining factor, $RC^0(D)$, is isomorphic to $RC(D)/\<\x_{S_0}\>$ and has  generators in one-to-one correspondence with those of $\pi_{dehn}(D)$. By Lemma \ref{flem} and Proposition \ref{freeproduct}, up to isomorphism $RC^0(D)$ does not depend on the particular region $S_0$.  
 
A fundamental region $\De$ of a genus-$g$ surface ${\cal S}$ has $4g$ boundary arcs, identified in pairs in the usual way. Assume that the diagram $D$ is oriented. Reading along an arc $\b$ in either direction, let $w_\b = x_{a_{i_1}}^{\e_1}\ldots x_{a_{i_r}}^{\e_r} \in \pi_{wirt}(D)$, where $a_{i_1}, \ldots, a_{i_r}$ are  the arcs of $D$ that are successively encountered, and $\e_i = a_i \cdot \b$  is the oriented intersection number computed in ${\cal S}$. Denote by $W(\De)$ the normal subgroup of $\pi_{wirt}(D)$ generated by the set of elements $w_\b$, with $\b$ ranging over the boundary arcs of $\De$.  By Theorem 3.2.2 of \cite{By12}   (also Theorem 4.4 of \cite{SW21}), $\pi_{wirt}(D)/W(\De)$ is isomorphic to $\pi_{dehn}(D)$. The isomorphism in \cite{SW21} identifies Wirtinger generators with the arc elements or their inverses in the Dehn group. 

\begin{remark} Reading along $\b$ in an opposite direction produces the inverse element. Hence the direction that we choose is unimportant. Also, if $\b'$ and $\b$ are identified in the surface, then $\omega_\b$ is equal to $\omega_{\b'}$ or $\omega^{-1}_{\b'}$. Therefore, we need consider only $2g$ boundary arcs. \end{remark} 

Let $AC(\De)$ be the normal subgroup of $AC(D)$ generated by the elements $\omega_\b = g_{a_{i_1}}g_{a_{i_2}}^{-1} g_{a_{i_3}}\ldots g_{a_{i_r}}^{-1}$ and $\omega'_\b= g_{a_{i_1}}^{-1}g_{a_{i_2}} g_{a_{i_3}}^{-1}\ldots g_{a_{i_r}}$, with $\b$ ranging over the boundary arcs
of $\De$.   \medskip

The following generalizes Theorem \ref{freeproduct}.
 
 \begin{theorem} \label{2.5+} If $D$ is a checkerboard colorable link diagram in a  surface ${\cal S}$ with fundamental region $\De$, then $AC(D)/AC(\De) \cong RC^0(D)$.
 \end{theorem} 
 
 \begin{proof} Recall from Section 2 that we have a function  $A(D)$ to $RC(D)$ sending $a$ to $\g_U^{}  \g_S^{-1}$, where $U$ (resp. $S$) is the unshaded 
 (resp. shaded) region incident at some point of the arc $a$. Moreover, the function determines a surjection $f: AC(D) \to RC^0(D)$. It is clear that $f$ vanishes on the subgroup $AC(\De)$, and hence induces a homomorphism $ \bar f: AC(D)/AC(\De)\to RC^0(D)$. 
 
 We have also a map $h: RC^0(D) \to AC(D)/AC(\De)$ sending the generator $\g_R$ of any region of $D$ to an element of $AC(D)$. We review the process. Recall we defined $\chi_{S_0} \in \pi_{dehn}(D)$ to be trivial in the definition of the Dehn group. Choose  a path from $S_0$ to $R$, crossing the arcs of $D$ transversely and avoiding crossings. Suppose that $a_0, a_1, \ldots, a_m$ are the arcs that we cross successively. If $R$ is shaded, then $h(\g_R)= g_{a_m}^{-1} g_{a_{m-1}} g_{a_{m-2}}^{-1} \cdots g_{a_0}$. If $R$ is unshaded then $h(\g_R)=g_{a_m} g_{a_{m-1}}^{-1} g_{a_{m-2}} \cdots g_{a_0}$.   
 
The process we have described is a sort of integration, accumulating generators of $AC(D)$ or their inverse each time we cross an arc. 
We can apply it to a path from one region, with an initial value, to another. We will refer to the final value that we obtain as the ``return value."  For example, if the path is a loop around a crossing of $D$, then one checks easily that the return value is equal to the initial value.  We can see why $h(\g_R)$ is independent of the choice of path from $S_0$ to $R$,  modulo $AC(\De)$, by considering an arbitrary closed path beginning and ending in $S_0$. Such a loop can be moved past any crossing without changing its return value. We can deform the loop to the boundary of $\De$, an embedded bouquet of loops in ${\cal S}$. 
Reading as above along any closed boundary arc returns the initial value modulo the subgroup $AC(\De)$. Hence the closed path with trivial initial value also returns the trivial value in the quotient $AC(D)/AC(\De)$. 
 
 Finally, one checks as in Section 2 that the maps $f$ and $h$ are inverses of each other.  \end{proof} 
 
\begin{remark}\label{basing} The definition of the map $h$ depends on the region $S_0$ for which $\chi_{S_0}$ is the identity. For the following, we will choose $S_0$ also to be the region for which $g_{S_0}$ is the identity in the definition of the Dehn group $\pi_{dehn}(D)$. \end{remark}
 
\begin{theorem} \label{catrcore} If $D$ is a checkerboard colorable diagram in a surface, then $\C(\pi_{dehn}(D)) \cong RC^0(D)$. \end{theorem}

\begin{proof}  Provide an orientation for $D$. The standard Wirtinger 3-tuple $(\pi_{wirt}(D), S, \rho)$ induces a  3-tuple $(\pi_{wirt}(D)/W(\De), \bar S, \bar \rho)$, with $\bar S$ and $\bar \rho$ defined naturally. Since an isomorphism between $\pi_{wirt}(D)/W(\De)$ and $\pi_{dehn}(D)$ matches Wirtinger generators with arc elements of $\pi_{dehn}(D)$ or their inverses, $(\pi_{wirt}(D)/W(\De), \bar S, \bar \rho)$  is isomorphic in the category {\bf Core} to the standard Dehn 3-tuple. Hence it suffices to prove that $\C(\pi_{wirt}(D)/W(\De), \bar S, \bar \rho)$ is isomorphic to $RC^0(D)$. 

The construction in the proof of Theorem \ref{catcore} yields a presentation of $\C(\pi_{wirt}(D)/W(\De), \bar S, \bar \rho)$. 
The cover $\widehat X$ that we produce with fundamental group $\C(\pi_{wirt}(D))$ acquires extra 2-cells corresponding to the generators of $W(\De)$. 
By Theorem \ref{catcore}, $\C(\pi_{wirt}(D)) \cong AC(D)$, and the boundaries of the extra 2-cells generate $AC(\De)$. 
Hence $\C(\pi_{wirt}(D)/W(\De), \bar N, \bar \rho)) \cong AC(D)/AC(\De)$. By Theorem \ref{2.5+}, the latter group is $RC^0(D).$ 
\end{proof} 

\begin{example} The diagram $D$ below represents a knot in a thickened torus. The Wirtinger group $\pi_{wirt}(D)$ has presentation 
\begin{equation}\label{wirt} \< a, b, c ; a b^{-1}a= b a b^{-1}, d^{-1} bd=b^{-1}d b, b^{-1}d b a^{-1}b^{-1} d b a b^{-1} d^{-1}ba^{-1} \>. \end{equation}
Adding the relations $a=d$ and $b=c$, read along the boundary of the fundamental region, produces a presentation for the Dehn group $\pi_{dehn}(D)$: 
\begin{equation} \label{dehn} \<a, b ; aba=bab, a^2 = b^2\>. \end{equation}
A presentation for $\pi_{dehn}(D)$ can also be calculated directly from the diagram: 
\begin{equation}\label{dehndirect} \<A, B ; ABA=BAB, A^2 = B^2\>.\end{equation} 
The group presented can be seen to be the semidirect product of $\Z$ acting nontrivially on $\Z_3$. 
Applying Proposition \ref{compute} to the presentation (\ref{wirt}) and using the substitutions $U = bd^{-1}, V= ba^{-1}$, we obtain the following presentation for $AC(D)$:
$$\< b, U, V ; U^3, V^3, (UV)^3\> \cong \Z * \<U, V; U^3, V^3, (UV)^3\>.$$
The second free factor is the plane crystallographic (``wallpaper") group {\bf p3} (see for example \cite{CM}, page 48). 
We apply Theorem \ref{2.5+} and Proposition \ref{compute} to the presentation (\ref{dehn}) to obtain a presentation of $RC^0(D)$:
$$\<a, b; (ab^{-1})^3 \> \cong \Z*\Z_3.$$
Hence $RC(D) \cong \Z*\Z*\Z_3$. 
Alternatively, one can compute a presentation of $RC(D)$ directly from presentation (\ref{dehndirect}): 
$$\<A, B, C; (BA^{-1})^3\>.$$

\begin{figure}[H]
\begin{center}
\includegraphics[height=2.5 in]{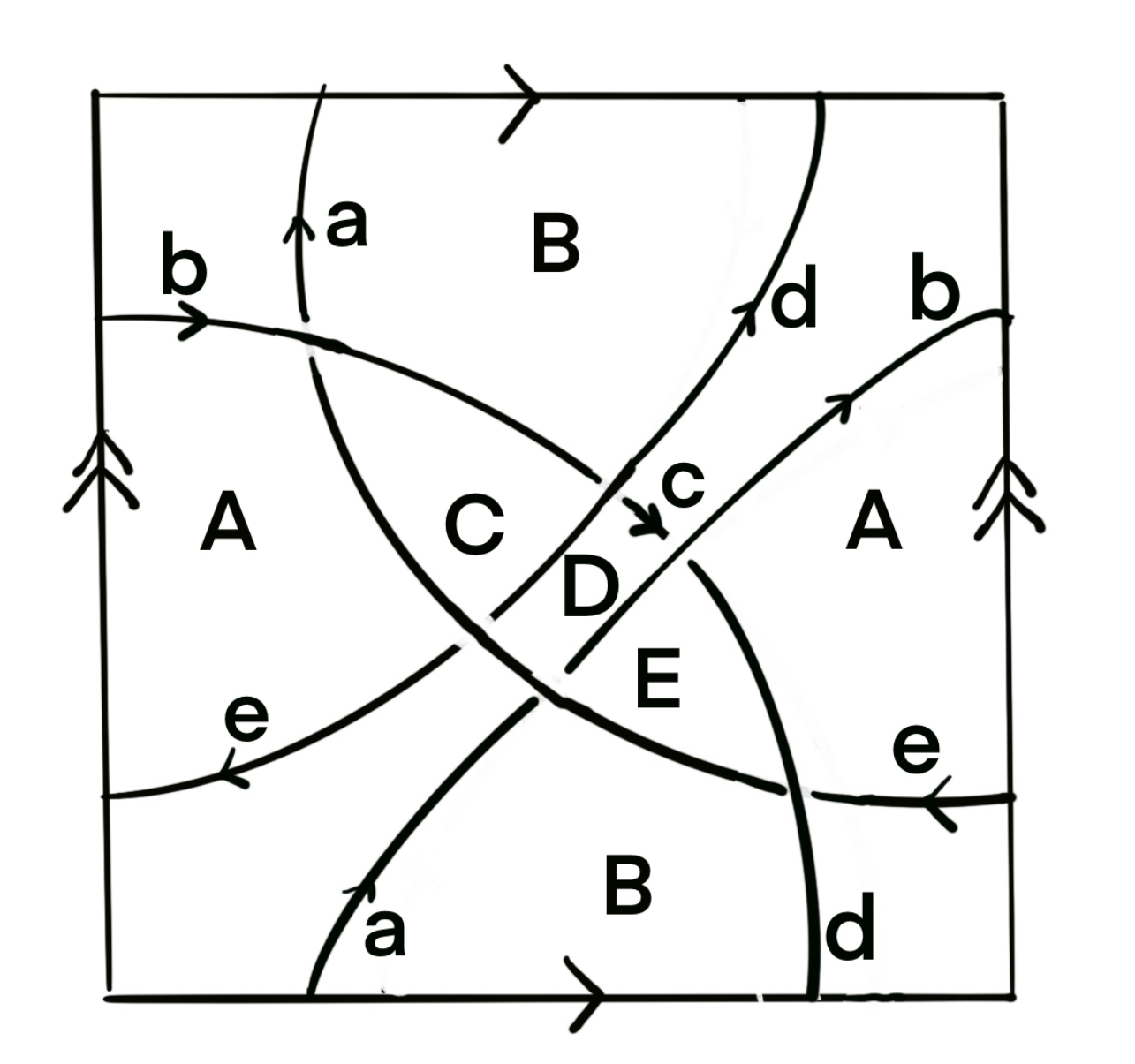}
\caption{$AC(D) \cong \Z *{\bf p3}$ and $RC(D) \cong \Z*\Z*\Z_3$}
\end{center}
\end{figure}

\end{example} 

We finish with the observation that the conclusion of Theorem \ref{catcore} does not hold in general for diagrams that are not checkerboard colorable. One can easily verify that for the diagram in the figure below the Dehn group $\pi_{dehn}(D)$ is trivial and hence so is $\C(\pi_{dehn}(D))$. However,  $RC^0(D)$ is cyclic of order 3.  

\begin{figure}[H]
\begin{center}
\includegraphics[height=2.2 in]{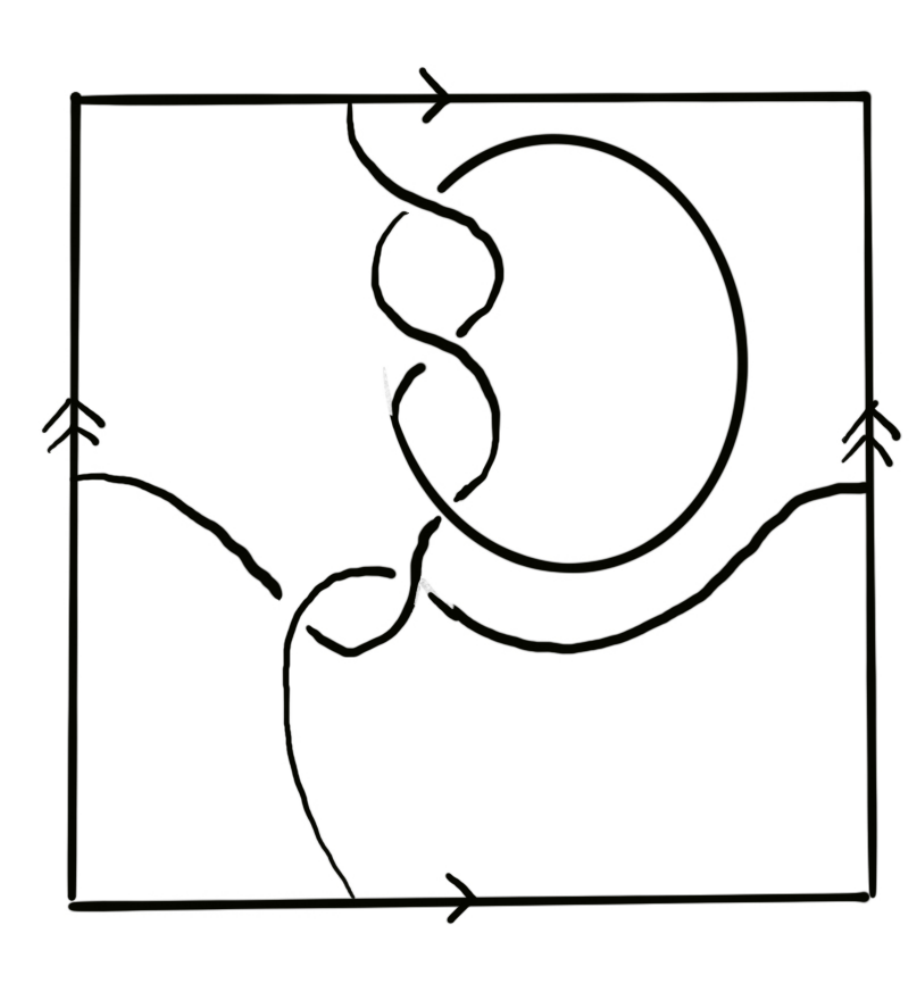}
\caption{Non-checkerboard colorable diagram $D$ such that $\C(\pi_{dehn}(D))$ and $RC^0(D)$ are distinct.}
\end{center}
\end{figure}

\bigskip

\ni Department of Mathematics\\
\ni Lafayette College\\Easton PA 18042\\
\ni Email: traldil@lafayette.edu\\ 

\ni Department of Mathematics and Statistics,\\
\ni University of South Alabama\\ Mobile, AL 36688 USA\\
\ni Email: silver@southalabama.edu\\ swilliam@southalabama.edu
\end{document}